\documentclass[11pt]{amsart}
\usepackage{amsfonts}
\usepackage{amssymb,latexsym}
\usepackage{enumerate}
\usepackage{mathrsfs}
\usepackage{amsmath}
\usepackage{amsthm}

scaled\magstephalf

\newtheorem{Thm}{Theorem}[section]

\newtheorem{Lem}[Thm]{Lemma}
\newtheorem{Rem}[Thm]{Remark}

\numberwithin{equation}{section}

\begin{document}
\setlength{\baselineskip}{1.2\baselineskip}
\title[Gradient Estimates]{A New Proof of Gradient Estimates for  Mean Curvature Equations with  Oblique Boundary Conditions}
\author{Jinju Xu}
\address{Department of Mathematics\\
         University of Science and Technology of China\\
         Hefei Anhui 230026 CHINA}
         \email{july25@mail.ustc.edu.cn}


\maketitle

\begin{abstract}
In this paper, we will use the maximum principle to give  a new proof of  the gradient estimates for mean curvature equations with some oblique derivative problems. Specially, we shall give a new proof for the capillary problem with zero gravity in any  dimension $n\geq 2$ and  Neumann problem  in $n=2, 3$ dimensions.
\end{abstract}

\section{Introduction}
 The interior gradient estimates and
 the Dirichlet problem for the prescribed mean curvature equation  have been extensively studied, see Gilbarg and Trudinger \cite{GT01}.

 Many authors have also considered  various oblique boundary value problems for  second order elliptic equations.
We refer to the  literature  Lieberman \cite{Lieb13} and the references  therein.

In this note, we mainly consider
 the following  oblique boundary value  problem for prescribed mean curvature equation
 \begin{align}
 \texttt{div}(\frac{Du}{\sqrt{1+|Du|^2}}) =&f(x, u)   \quad\text{in}\quad \Omega, \label{1.1}\\
              v^{q-1}\frac{\partial u}{\partial \gamma}+\psi(x, u)=&0  \quad\text{on} \quad\partial \Omega,\label{1.2}
\end{align}
where $\Omega\subset\mathbb R^n$ is a bounded domain, $n\geq 2$, $\gamma$ is the inward unit normal to $\partial\Omega$ and $q\geq 0$, $v=\sqrt{1+|Du|^2}$.

In \eqref{1.2},  for $q=0$,  it is corresponding to capillary  boundary condition and for $q=1$,  it is corresponding to Neumann boundary value.

   For the  mean curvature equation with capillary problem, there have been many existence results such as Ural'tseva \cite{Ur73}, Simon-Spruck \cite{SS76}, Gerhardt \cite{Ger76}. They obtained gradient  estimates  via test function technique. Spruck\cite{Sp75} used the maximum principle to obtain boundary gradient estimate in two dimension for positive gravity case ( $f_u\ge C_0>0$, $C_0$ is a constant). Korevaar\cite{Kor88} generalized his normal variation technique and got the gradient estimates in the positive gravity case in high dimensions case.   Simultaneously,  Lieberman\cite{Lieb88} used maximum principle to get  the gradient estimates on  general quasilinear elliptic equations with capillary problem in zero gravity case ($f_u\ge 0$).

  For the problem \eqref{1.1}, \eqref{1.2},  Lieberman (\cite{Lieb13}, in page 360) proved the gradient estimates  for $q>1$  or $q=0$. Recently, Ma-Xu\cite{MX14} have given the  gradient estimates of  mean curvature equations with Neumann problem via maximum principle. Moreover, they got an existence result in positive gravity case.

In this paper, we use maximum principle  to give new proofs of gradient estimates for  the problem \eqref{1.1}, \eqref{1.2} with $q>1$ or $q=0$ or $q=1$ ( $n=2,3$ dimensions ) cases respectively. Our proofs are elementary and  based on the choice of auxiliary functions.

Let's restate the following three results. First consider the boundary value condition with $\psi=\psi(x)$.
\begin{align}
v^{q-1}\frac{\partial u}{\partial \gamma}+\psi(x)=0  \quad\text{on} \quad\partial \Omega.\label{1.2a}
\end{align}
\begin{Thm}[\cite{Lieb13}]\label{Thm1.1}  Let  $\Omega\subset\mathbb R^n$ be a bounded domain, $n\geq 2$, $\partial \Omega \in C^{3}$, and $\gamma$ be the inward unit normal to $\partial\Omega$.
Suppose $u\in C^{2}(\overline\Omega)\bigcap C^{3}(\Omega)$ is a  solution of \eqref{1.1}, \eqref{1.2a} with $q>1$.
$f(x,z), \psi(x)$ are given functions defined in $\overline\Omega\times [-M_0, M_0]$ and $ \overline\Omega$ respectively. Furthermore we  assume
there exist positive constants $M_0, L_1, L_2$ such that
\begin{align}
|u|\leq& M_0\quad \text{in}\quad\overline\Omega,\label{1.3}\\
f_z(x,z)\geq &0 \quad \text{in}\quad\overline\Omega\times[-M_0, M_0],\label{1.4}\\
|f(x,z)|+|f_{x}(x,z)|\leq & L_1 \quad \text{in}\quad \overline\Omega\times[-M_0, M_0],\label{1.5}\\
|\psi(x)|_{C^{1}(\overline\Omega)}\leq&  L_2.\label{1.6}
\end{align}
Then there exists a small positive constant
$\mu_0$ such that
$$\sup_{\overline\Omega_{\mu_0}}|Du|\leq \max\{M_1, M_2\},$$
where $M_1$ is a positive constant depending only on $n, \mu_0, M_0, L_1$, which is from the interior gradient estimates;
$M_2$ is  a positive constant depending only on $n, \Omega, \mu_0, M_0, L_1, L_2$, and $d(x) =\texttt{dist}(x, \partial\Omega), \Omega_{\mu_0} = \{x \in \Omega: d(x)<\mu_0\}.$
\end{Thm}

\begin{Thm}[\cite{Ger76}, \cite{Lieb88}]\label{Thm1.2} Let  $\Omega\subset\mathbb R^n$ be a bounded domain, $n\geq 2$, $\partial \Omega \in C^{3}$, and $\gamma$ be the inward unit normal to $\partial\Omega$.
Suppose $u\in C^{2}(\overline\Omega)\bigcap C^{3}(\Omega)$ is a  solution of \eqref{1.1}, \eqref{1.2a} with $q=0$ and satisfies \eqref{1.3}.
$f(x,z), \psi(x)$ are given functions defined in $\overline\Omega\times [-M_0, M_0]$ and $ \overline\Omega$ respectively. Assume $f(x,z)$ satisfies \eqref{1.4}-\eqref{1.5} and $\psi(x)$ satisfies \eqref{1.6}. Furthermore we assume there exists a positive constant $b_0$ such that
\begin{align}
 |\psi(x)|_{C^{0}(\partial\Omega)}\leq&  b_0<1.\label{1.7}
\end{align}
Then there exists a small positive constant
$\mu_0$ such that
$$\sup_{\overline\Omega_{\mu_0}}|Du|\leq \max\{M_1, M_2\},$$
where $M_1$ is a positive constant depending only on $n, \mu_0, M_0, L_1$, which is from the interior gradient estimates;
$M_2$ is  a positive constant depending only on $n, \Omega, \mu_0, M_0, L_1, L_2, b_0$.
\end{Thm}

The  following boundary gradient estimate of solutions for Neumann problem of mean curvature equations has been given by Ma-Xu\cite{MX14} in any dimension $n\geq 2$.
\begin{Thm}[\cite{MX14}]\label{Thm1.3} Let $\Omega\subset\mathbb R^n$ (n=2,3) be a bounded domain,  $\partial \Omega \in C^{3}$, and $\gamma$ be the inward unit normal to $\partial\Omega$.
Suppose $u\in C^{2}(\overline\Omega)\bigcap C^{3}(\Omega)$ is a  solution of \eqref{1.1}, \eqref{1.2} with $q=1$  and satisfies \eqref{1.3}.
$f(x,z), \psi(x,z)$ are given functions defined in $(\overline\Omega\times [-M_0, M_0])$ respectively. And $f(x,z)$ satisfies the conditions \eqref{1.4}-\eqref{1.5}. Furthermore assume
there exists a positive constant $L_3$ such that
\begin{align}
|\psi(x,z)|_{C^{1}(\overline\Omega\times [-M_0, M_0])}\leq&  L_3.\label{1.8}
\end{align}
Then
there exists a small positive constant
$\mu_0$ such that we have the following estimate
$$\sup_{\overline\Omega_{\mu_0}}|Du|\leq \max\{M_1, M_2\},$$
where $M_1$ is a positive constant depending only on $n, \mu_0, M_0, L_1$, which is from the interior gradient estimates;
$M_2$ is  a positive constant depending only on $n, \Omega, \mu_0, M_0, L_1, L_3$.
\end{Thm}

 As we stated before, there is a standard interior gradient estimates for the mean curvature equation.

\begin{Rem}[\cite{GT01}]\label{Rem1.1}
If $u\in C^{3}(\Omega)$ is a bounded solution for the equation \eqref{1.1}  with \eqref{1.3}, and if $f \in C^{1}(\overline\Omega \times [-M_0, M_0])$ satisfies the conditions \eqref{1.4}-\eqref{1.5}, then for any subdomain
$\Omega'\subset\subset\Omega$, we have
$$\sup_{\Omega'}|Du|\leq M_1,$$
where $M_1$ is a positive constant depending only on $n,  M_0, \texttt{dist} (\Omega', \partial\Omega), L_1$.
\end{Rem}

The rest of the paper is organized as follows. In section 2, we first give the definitions and some notations. We prove  Theorem~\ref{Thm1.1} in section 3 under the help of one lemma. This lemma  will be proved in section 4. In section 5, we give the proof of  Theorem~\ref{Thm1.2}. Finally we give the proof of Theorem~\ref{Thm1.3} in section 6.

\section{PRELIMINARIES}
We denote by $\Omega$ a bounded  domain in $\mathbb{R}^n$, $n\geq 2$,  $\partial \Omega\in C^{3}$,   set
\begin{align*}
 d(x)=\texttt{dist}(x,\partial \Omega),
 \end{align*}
 and
\begin{align*}
 \Omega_\mu=&\{{x\in\Omega:d(x)<\mu}\}.
 \end{align*}
Then it is well known that there exists a positive constant $\mu_{1}>0$ such that $d(x) \in C^3(\overline \Omega_{\mu_{1}})$. As in Simon-Spruck \cite{SS76} or Lieberman \cite{Lieb13} in page 331,  we can take $\gamma= D d$ in $\Omega_{\mu_{1}}$ and note that  $\gamma$ is a $C^2(\overline \Omega_{\mu_{1}})$ vector field. As mentioned in \cite{Lieb88} and the book \cite{Lieb13}, we also have the following formulas

\begin{align}\label{2.1}
\begin{split}
|D\gamma|+|D^2\gamma|\leq& C(n,\Omega) \quad\text{in}\quad \Omega_{\mu_{1}},\\
 \sum_{1\leq i\leq n}\gamma^iD_j\gamma^i=0,  \sum_{1\leq i\leq n}\gamma^iD_i\gamma^j=&0, \,|\gamma|=1 \quad\text{in} \quad\Omega_{\mu_{1}}.
\end{split}
\end{align}
As in \cite{Lieb13}, we define
 \begin{align}\label{2.2}
\begin{split}
c^{ij}=&\delta_{ij}-\gamma^i\gamma^j  \quad \text{in} \quad \Omega_{\mu_{1}},
\end{split}
\end{align}
 and for a vector $\zeta \in R^n$, we write $\zeta'$ for the vector with $i-$th component $ \sum_{1\leq j\leq n}c^{ij}\zeta_j$. So
 \begin{align}\label{2.3}
\begin{split}
|D'u|^2=& \sum_{1\leq i,j\leq n}c^{ij}u_iu_j.
\end{split}
\end{align}
Let
\begin{align}\label{2.4}
\begin{split}
a^{ij}(D u)=v^2\delta_{ij}-u_iu_j, \quad
v=(1+|D u|^2)^{\frac{1}{2}}.
\end{split}
\end{align}

Then the equations \eqref{1.1}, \eqref{1.2}  are equivalent to the following boundary value problem
\begin{align}
\sum_{i,j=1}^n a^{ij}u_{ij}=&f(x,u)v^3 \quad \text{in}\quad \Omega,\label{2.5}\\
u_{\gamma}=&-v^{1-q} \psi(x,u) \quad \text{on}\quad\partial\Omega.\label{2.6}
\end{align}

\section{Proof of Theorem~\ref{Thm1.1}. }

Now we begin to prove Theorem~\ref{Thm1.1}.,
 using the technique developed by Spruck \cite{Sp75}, Lieberman \cite{Lieb88} and Wang \cite{Wang98},
 we shall choose an auxiliary function which contains  $|D'u|^2$ and other lower order terms. Then we use the maximum principle for this auxiliary function in $\overline\Omega_{\mu_0}, 0<\mu_0<\mu_1$. At last,  we get our estimates.

{\em Proof of  Theorem~\ref{Thm1.1}.}
Let
$$P(x)=\log|D'u|^2e^{\sqrt{n}\alpha_0(M_0+1+u)}e^{\alpha_0d},$$ where we  have let
$$\alpha_0=2C_0L_2+2C_0+2, $$ which is a constant, and
$C_0$ is also a positive constant depending only on $n,\Omega$.

In order to simplify the computation, let
$$\varphi(x)=\log P(x)=\log\log|D'u|^2+h(u)+g(d),$$
where in the $q>1$ boundary value case, we choose
\begin{align}\label{3hg1}
h(u)=\sqrt{n}\alpha_0(M_0+1+u), \quad g(d)=\alpha_0 d.
\end{align}

 We assume that
$\varphi(x)$ attains its maximum at $x_0 \in \overline \Omega_{\mu_{0}}$, where $0<\mu_0<\mu_1$ is a sufficiently small number which we shall decide it  later.

Now we divide three cases to complete the proof of  Theorem~\ref{Thm1.1}.

Case I. If $\varphi(x)$ attains its maximum at $x_0 \in \partial\Omega$, then we shall use the Hopf Lemma to get the bound of $|D'u|(x_0)$.

Case II. If $\varphi(x)$ attains its maximum at $x_0 \in\partial\Omega_{\mu_0}\bigcap\Omega$, then we shall get the estimates via the standard interior gradient bound \cite{GT01}.

Case III.  If $\varphi(x)$ attains its maximum at $x_0 \in \Omega_{\mu_0}$, in this case for the sufficiently small constant $\mu_0>0$,  then we can use the maximum principle to get the bound of $|D'u|(x_0)$.

Now  all computations work at the point $x_0$.

{\bf Case I.} If $x_0\in \partial \Omega$, we shall get the bound of $|D'u|(x_0)$.

We differentiate $\varphi$ along the normal direction.
\begin{align}\label{3.1}
\frac{\partial\varphi}{\partial\gamma}=&\frac{\sum_{1\leq i\leq n}(|D'u|^2)_i\gamma^i}{|D'u|^2\log|D'u|^2}+h'u_{\gamma}+g'.
\end{align}
Applying \eqref{2.1} and \eqref{2.3}, it follows that
\begin{align}\label{3.2}
\sum_{1\leq i\leq n}(|D'u|^2)_i\gamma^i=&\sum_{1\leq i\leq n}(\sum_{1\leq k,l\leq n}c^{kl}u_{k}u_l)_i\gamma^i
=2\sum_{1\leq i,k,l\leq n}c^{kl}u_{ki}u_l\gamma^i.
\end{align}
Differentiating \eqref{2.6} with respect to tangential direction,   we have
\begin{align}\label{3.3}
\sum_{1\leq k\leq n}c^{kl}(u_{\gamma})_k=&-\sum_{1\leq k\leq n}c^{kl}(v^{1-q}\psi )_k.
\end{align}
It follows that
\begin{align}\label{3.4}
\sum_{1\leq i, k\leq n}c^{kl}u_{ik}\gamma^i=&-\sum_{1\leq i,k\leq n}c^{kl}u_i(\gamma^i)_k-v^{1-q}\sum_{1\leq k\leq n}c^{kl}D_k\psi-(1-q)\psi v^{-q}\sum_{1\leq k\leq n}c^{kl}v_k.
\end{align}
Here in order to avoid repeated calculation in the back, we have let $\psi= \psi(x,u)$ and
\begin{align}\label{3Dkpsi}
D_k\psi=&\psi_{x_k}+\psi_uu_k.
\end{align}
Since
\begin{align}\label{3.5}
v^2=&1+|D'u|^2+u^2_\gamma,
\end{align}
differentiating \eqref{3.5} with respect to $x_k$, we obtain
\begin{align}\label{3.6}
v_k=&\frac{(|D'u|^2)_k}{2v}+\frac{u_{\gamma}}{v}\sum_{1\leq i\leq n}(u_{ik}\gamma^i+u_i(\gamma^i)_k).
\end{align}
Inserting \eqref{3.6} into \eqref{3.1}, we have
\begin{align}\label{3.7}
\begin{split}
\sum_{1\leq i,k\leq n}c^{kl}u_{ki}\gamma^i=&-\sum_{1\leq i,k\leq n}c^{kl}u_i(\gamma^i)_k-\frac{v^{1-q}\sum_{1\leq k\leq n}c^{kl}D_k\psi}{1-(1-q)v^{-2q} \psi^2}\\
&-\frac{1-q}{2}\cdot\frac{\psi v^{-q-1}}{1-(1-q)v^{-2q} \psi^2} \sum_{1\leq k\leq n}c^{kl}(|D'u|^2)_k.
\end{split}
\end{align}
Since $\sum_{1\leq k\leq n}c^{kl}\varphi_k=0$,  and $\sum_{1\leq k\leq n}c^{kl}\gamma^k=0$, we obtain
\begin{align}\label{3.8}\begin{split}
\sum_{1\leq k\leq n}c^{kl}(|D'u|^2)_k=&-|D'u|^2\log|D'u|^2\sum_{1\leq k\leq n}c^{kl}(h'u_k+g'\gamma^k)\\
=&-h'|D'u|^2\log|D'u|^2\sum_{1\leq k\leq n}c^{kl}u_k.
\end{split}\end{align}
Inserting \eqref{3.8} into \eqref{3.7}, we have
\begin{align}\label{3.9}
\begin{split}
\sum_{1\leq i,k\leq n}c^{kl}u_{ki}\gamma^i=&-\sum_{1\leq i,k\leq n}c^{kl}u_i(\gamma^i)_k-\frac{v^{1-q}\sum_{1\leq k\leq n}c^{kl}D_k\psi}{1-(1-q)v^{-2q} \psi^2}\\
&+\frac{1-q}{2}\cdot\frac{h'\psi v^{-q-1}}{1-(1-q)v^{-2q} \psi^2}|D'u|^2\log |D'u|^2 \sum_{1\leq k\leq n}c^{kl}u_k.
\end{split}
\end{align}
Putting \eqref{3.9} into \eqref{3.2}, combining \eqref{3.1}, we have
\begin{align}\label{3.10}
\begin{split}
|D'u|^2\log |D'u|^2 \frac{\partial\varphi}{\partial\gamma}
=&g'|D'u|^2\log |D'u|^2-2\sum_{1\leq i,k,l\leq n}c^{kl}u_iu_l(\gamma^i)_k-\frac{2v^{1-q}\sum_{1\leq k,l\leq n}c^{kl}D_k\psi u_l}{1-(1-q)v^{-2q} \psi^2}\\&
-h'\psi \frac{q v^{1-q}+(1-q)v^{-1-q}}{1-(1-q)v^{-2q} \psi^2}|D'u|^2\log |D'u|^2.
\end{split}
\end{align}
In the following, we consider $\psi=\psi(x)$, $q>1$, then we have
\begin{align}\label{3.11}
\begin{split}
|D'u|^2\log |D'u|^2 \frac{\partial\varphi}{\partial\gamma}
=&g'|D'u|^2\log |D'u|^2-2\sum_{1\leq i,k,l\leq n}c^{kl}u_iu_l(\gamma^i)_k-\frac{2v^{1-q}\sum_{1\leq k,l\leq n}c^{kl}\psi_k u_l}{1-(1-q)v^{-2q} \psi^2}\\&
-h'\psi \frac{q v^{1-q}+(1-q)v^{-1-q}}{1-(1-q)v^{-2q} \psi^2}|D'u|^2\log |D'u|^2.
\end{split}
\end{align}
Since at $x_0$ , $$\psi^2=v^{2q-2}u_\gamma^2=(1+|Du|^2)^{q-1}u_\gamma^2, \quad |Du|^2=|D'u|^2+u_\gamma^2, $$
then we have
if
\begin{align}\label{3.12a}
b_1|D'u|^2<u_\gamma^2,
\end{align}
where $b_1$ is a positive constant,
then \begin{align}\label{3.12b}
u_\gamma^2>\frac{b_1}{1+b_1}|Du|^2,
\end{align}
and we get the estimate
\begin{align}\label{3.12c}
(1+|Du|^2)^{q-1}|Du|^2<\frac{1+b_1}{b_1}\psi^2,
\quad |Du|^2<(\frac{1+b_1}{b_1}\psi^2)^{\frac{1}{q}},
\end{align}
 and we complete this proof.

So we can assume
\begin{align}\label{3.12d}
b_1|D'u|^2\geq u_\gamma^2,
\end{align}
 then from  $|Du|^2=|D'u|^2+u_\gamma^2$, we have
 \begin{align}\label{3.12e}|Du|^2\leq (1+b_1)|D'u|^2.
 \end{align}
Now we assume at $x_0$,  we have
\begin{align}\label{3.12f}
|Du|\geq \max\{ 10\sqrt{1+b_1}, (4q\sqrt{n}|\psi|_{C^{0}(\partial\Omega)})^{\frac{1}{q-1}} \},
 \end{align}
 then we can get the the following estimate at $x_0$,
 \begin{align}\label{3.12g}
 |D'u|\geq \max\{10, \frac{(4q\sqrt{n}|\psi|_{C^{0}(\partial\Omega)})^{\frac{1}{q-1}}}{\sqrt{1+b_1}}\}.
 \end{align}
Inserting  \eqref{3.12g} into \eqref{3.11}, and by the choice of $h(u), g(d)$ in \eqref{3hg1},  we obtain
\begin{align}\label{3.13}
\begin{split}
|D'u|^2\log |D'u|^2 \frac{\partial\varphi}{\partial\gamma}
\geq&\big[\alpha_0-C_0-C_0|\psi|_{C^{1}(\overline{\Omega})}
-\frac{2q\sqrt{n}|\psi|_{C^{0}(\partial\Omega)}}{v^{q-1}}\alpha_0\big]|D'u|^2\log |D'u|^2\\
\geq&|D'u|^2\log |D'u|^2\\
>&0.
\end{split}
\end{align}
On the other hand, by the Hopf Lemma, we have
$$\frac{\partial\varphi}{\partial\gamma}(x_0)\leq 0,$$ it is a contradiction to \eqref{3.13}.
Then we have
$$|D'u|(x_0)\leq\max\{10, \frac{(4q\sqrt{n}|\psi|_{C^{0}(\partial\Omega)})^{\frac{1}{q-1}}}{\sqrt{1+b_1}}\}.$$
So
$$|Du|(x_0)\leq\max\{(\frac{1+b_1}{b_1}\psi^2)^{\frac{1}{2q}},\, 10\sqrt{1+b_1}, \, (4q\sqrt{n}|\psi|_{C^{0}(\partial\Omega)})^{\frac{1}{q-1}}\}.$$

{\bf Case II.} $ x_0\in \partial\Omega_{\mu_{0}}\bigcap\Omega$. This is due to interior gradient estimates. From Remark~\ref{Rem1.1}, we have
 \begin{align}\label{3CaseII}
\sup_{\partial\Omega_{\mu_0}\bigcap\Omega}|Du|\leq \tilde{M}_1.
\end{align}
where $\tilde{M}_1$ is a positive constant depending only on $n, M_0, \mu_0, L_1$.

{\bf Case III.} $x_0\in\Omega_{\mu_{0}}$. we shall get the bound of $|D'u|(x_0)$.\par

In this case, $x_0$ is a critical point of $\varphi$. We choose the normal coordinate at $x_0$, by rotating the coordinate system
suitably, we may assume that $u_i(x_0)=0,\,2\leq i\leq n$ and $u_1(x_0)=|Du|>0$. And we can further assume that the matrix $(u_{ij}(x_0))(2\leq i,j\leq n)$ is diagonal.

 We can choose $0<\mu_0< \mu_1$ , and $\mu_0$ is sufficiently small.

 From Gilbarg and Trudinger\cite{GT01} [page 368, formula (15.38)], we have
 \begin{align}\label{3.14a}
\sup_{\Omega}|Du|^2\leq C_1(1+\sup_{\partial\Omega}|Du|^2),
\end{align}
 where $C_1$ is a positive constant depending on $n, L_1, M_0$.

From  Case I,   we can assume \eqref{3.12d}, otherwise we have finished the proof of  Theorem~\ref{Thm1.1}.
Hence,
\begin{align}\label{3.14b}
\begin{split}
\sup_{\Omega}|Du|^2\leq& C_1\big[1+(1+b_1)\sup_{\partial\Omega\bigcap\{|D'u|\geq 1\}}|D'u|^2\big]\\
\leq& C_2\sup_{\partial\Omega\bigcap\{|D'u|\geq 1\}}|D'u|^2.
\end{split}
\end{align}

 So we have
\begin{align}\label{3.14c}
\begin{split}
\sup_{\Omega_{\mu_0}}|Du|^2
\leq& C_3\sup_{\overline{\Omega}_{\mu_0}(M)}|D'u|^2,
\end{split}
\end{align}
  where $\Omega_{\mu_0}(M)=\Omega_{\mu_0}\bigcap\{|D'u|\geq M\}$, $M>10$ is a positive constant;
   $C_3$ is a positive constant depending on $n, L_1,  M_0$.

Assume  $x_1\in\Omega_{\mu_0}(M)$ such that
\begin{align}\label{3.15a}
\begin{split}
\sup_{\Omega_{\mu_0}(M)}|D'u|^2
=& |D'u|^2(x_1).
\end{split}
\end{align}
Since $P(x_0)\geq P(x_1)$, then we have
\begin{align}\label{3.15b}
\begin{split}
\log|D'u|^2(x_0)h(u(x_0))g(d(x_0))\geq&\log|D'u|^2(x_1)h(u(x_1))g(d(x_1)).
\end{split}
\end{align}
It follows that
\begin{align}\label{3.15c}
\begin{split}
|D'u|^2(x_1)\leq&C_4|D'u|^2(x_0).
\end{split}
\end{align}
where  $C_4$ is a positive constant depending on $n, L_1, L_2,  M_0$.
However
\begin{align}\label{3.15d}
\begin{split}
\sup_{\Omega_{\mu_0}}|Du|^2
\leq& C_3\sup_{\Omega_{\mu_0}(M)}|D'u|^2=|D'u|^2(x_1)\leq C_4|D'u|^2(x_0).
\end{split}
\end{align}
 Assume $|D'u|(x_0)\geq M$, otherwise we get the estimate. Hence at $x_0$,
\begin{align}\label{3.15e}
\begin{split}
u_1^2(x_0)=|Du|^2(x_0)
\leq& C_4c^{11}u_1^2(x_0).
\end{split}
\end{align}
Then we have at $x_0$,
\begin{align}\label{3c^{11}}
\begin{split}
c^{11}
\geq& \frac{1}{C_4}>0.
\end{split}
\end{align}
From the above choices, we shall prove Theorem~\ref{Thm1.1} with three steps, as we mentioned before, all the calculations will be done at the fixed point $x_0$.

{\bf Step 1:} We first get the formula \eqref{3aijvarphiijb}.\par

Taking the first  derivative of $\varphi$,
\begin{align}\label{3varphii}
\varphi_i=&\frac{(|D'u|^2)_i}{|D'u|^2\log|D'u|^2}+h'u_i+g'\gamma^i.
\end{align}
From $\varphi_i(x_0)=0$,  we have
\begin{align}\label{3varphii=0}
(|D'u|^2)_i=-|D'u|^2\log|D'u|^2(h'u_i+g'\gamma^i).
\end{align}
Take the derivative again for  $\varphi_i$,
\begin{align}\label{3varphiija}
\begin{split}
\varphi_{ij}
=&\frac{(|D'u|^2)_{ij}}{|D'u|^2\log|D'u|^2}-(1+\log|D'u|^2)\frac{(|D'u|^2)_i(|D'u|^2)_j}{(|D'u|^2\log|D'u|^2)^2}\\&
+h'u_{ij}
+h''u_iu_j+g''\gamma^i\gamma^j+g'(\gamma^i)_j.
\end{split}
\end{align}
Using \eqref{3varphii=0}, it follows that
\begin{align}\label{3varphiijb}
\begin{split}
\varphi_{ij}
=&\frac{(|D'u|^2)_{ij}}{|D'u|^2\log|D'u|^2}+h'u_{ij}+\big[h''-(1+\log|D'u|^2)h'^2\big]u_iu_j\\&
+\big[g''-(1+\log|D'u|^2)g'^2\big]\gamma^i\gamma^j
-(1+\log|D'u|^2)h'g'(\gamma^iu_j+\gamma^ju_i)+g'(\gamma^i)_j.
\end{split}
\end{align}
Then we get
\begin{align}\label{3aijvarphiija}
\begin{split}
0\geq \sum_{1\leq i,j\leq n}a^{ij}\varphi_{ij}
=:&I_1+I_2,
\end{split}
\end{align}
where
\begin{align}\label{3I1a}
I_1=\frac{1}{|D'u|^2\log|D'u|^2}\sum_{1\leq i,j\leq n}a^{ij}(|D'u|^2)_{ij},
\end{align}
and
\begin{align}
I_2=&\sum_{1\leq i,j\leq n}a^{ij}\bigg\{h' u_{ij}+\big[h''-(1+\log|D'u|^2)h'^2\big]u_iu_j+\big[g''-(1+\log|D'u|^2)g'^2\big]\gamma^i\gamma^j\notag\\
&\qquad\qquad\quad-2(1+\log|D'u|^2)h'g'\gamma^iu_j+g'(\gamma^i)_j\bigg\}\label{3I2a}.
\end{align}
From the choice of the coordinate, we have
\begin{align}\label{3aii}
a^{11}=1,\, a^{ii}=v^2=1+u_1^2 \,\,(2\leq i\leq n ), a^{ij}=0\,\, (i\neq j,1\leq i,j\leq n).
\end{align}
and
\begin{align}\label{3|D'u|2a}
|D'u|^2=c^{11}u^2_1,\quad |D'u|^2\log|D'u|^2=2c^{11}u_1^2\log u_1+c^{11}(\log c^{11})u_1^2.
\end{align}

Now we first treat $I_2$.

From the equations \eqref{2.5}, \eqref{3aii}  and \eqref{3|D'u|2a}, we have
\begin{align}
I_2
=&h'f v^3-h'^2u_1^2\log|D'u|^2+(h''-h'^2)u_1^2+\big[g''-(1+\log|D'u|^2)g'^2\big]\sum_{1\leq i\leq n}a^{ii}(\gamma^i)^2\notag\\
&-2(1+\log|D'u|^2)h'g'\gamma^1u_1+g'\sum_{1\leq i\leq n}a^{ii}(\gamma^i)_i,\notag\\
=&h'f v^3-2(h'^2+c^{11}g'^2) u_1^2\log u_1+\big[h''-(1+\log c^{11})h'^2-c^{11}(1+\log c^{11})g'^2\notag\\
&+c^{11}g''+g'\sum_{2\leq i\leq n}(\gamma^i)_i\big] u_1^2-4h'g'\gamma^1 u_1\log u_1-2(1+\log c^{11})h'g'\gamma^1 u_1\notag\\
&-2g'^2 \log u_1+g''-(1+\log c^{11})g'^2
+g'\sum_{1\leq i\leq n}(\gamma^i)_i.\label{3I2b}
\end{align}
So we have
\begin{align}
I_2
\geq&\sqrt{n}\alpha_0f v^3-2(n+c^{11})\alpha_0^2 u_1^2\log u_1-C_5u_1^2,\label{3I2c}
\end{align}
here we use the expression for $h(u), g(d)$ in \eqref{3hg1}, and $C_5$  is a positive constant  depending  only on $n, \Omega, M_0, \mu_0, L_2$.\par

Next, we calculate $I_1$ and get the formula \eqref{3I1c}. \par

From \eqref{2.3}, taking the first  derivative of $|D'u|^2$, we have

\begin{align}\label{3|D'u|^2i}
\begin{split}
(|D'u|^2)_i=&\sum_{1\leq k, l\leq n}(c^{kl})_iu_ku_l+2\sum_{1\leq k, l\leq n}c^{kl}u_{ki}u_l.
\end{split}
\end{align}

Taking the derivatives of $|D'u|^2$ once more, we have

\begin{align}\label{3|D'u|^2ij}
\begin{split}
(|D'u|^2)_{ij}=&\sum_{1\leq k, l\leq n}(c^{kl})_{ij}u_ku_l+2\sum_{1\leq k, l\leq n}(c^{kl})_iu_{kj}u_l+2\sum_{1\leq k, l\leq n}(c^{kl})_ju_{ki}u_l\\
&+2\sum_{1\leq k, l\leq n}c^{kl}u_{kij}u_l+2\sum_{1\leq k, l\leq n}c^{kl}u_{ki}u_{lj}.
\end{split}
\end{align}
By \eqref{3I1a} and \eqref{3|D'u|^2ij}, we can rewrite $I_1$ as
\begin{align}\label{3I1b}
\begin{split}
I_1
=&\frac{1}{|D'u|^2\log|D'u|^2}\big[I_{11}+I_{12}+I_{13}+I_{14}\big],
\end{split}
\end{align}
where
\begin{align*}
I_{11}=&u^2_1\sum_{1\leq i\leq n}a^{ii}(c^{11})_{ii} ,\quad
I_{12}=4u_1\sum_{1\leq i,k\leq n}a^{ii}(c^{k1})_{i}u_{ki} ,\\
I_{13}=&2u_1\sum_{1\leq i,j,k\leq n}c^{k1}a^{ij}u_{ijk},\quad
I_{14}=2\sum_{1\leq i,k,l\leq n}c^{kl}a^{ii}u_{ki}u_{li}.
\end{align*}

In the following, we shall deal with $I_{11}, I_{12}, I_{13}$ and $I_{14}$ respectively. \par
For the terms $I_{11}$ and $I_{12}$: from  \eqref{3aii}, we have
\begin{align}\label{3I11}
I_{11}=&\sum_{2\leq i\leq n}(c^{11})_{ii} u^4_1+\sum_{1\leq i\leq n}(c^{11})_{ii}u^2_1,\\
I_{12}
=&4(c^{11})_{1} u_1u_{11}+4 u_1\sum_{2\leq i\leq n}[(c^{1i})_{1}+v^2(c^{11})_{i}]u_{1i}+4u_1v^2\sum_{2\leq i\leq n}(c^{1i})_{i}u_{ii}.\label{3I12}
\end{align}

For the term $I_{13}$:
by the equation \eqref{2.5}, we have
\begin{align}\label{3u11a}
\begin{split}
u_{11}=&fv^3-v^2\sum_{2\leq i\leq n}u_{ii},
\end{split}
\end{align}
and
\begin{align}\label{3Deltau}
\begin{split}
\Delta u=&fv+\frac{u_1^2}{v^2}u_{11}.
\end{split}
\end{align}
Differentiating \eqref{2.5}, we have
\begin{align}\label{3aijuijka}
\begin{split}
\sum_{1\leq i,j\leq n}a^{ij}u_{ijk}=&-\sum_{1\leq i,j,l\leq n}a^{ij}_{p_l}u_{lk}u_{ij}+v^3D_{k}f+3fv^2v_k.
\end{split}
\end{align}
From \eqref{2.4}, we have
\begin{align}\label{3aijpl}
\begin{split}
a^{ij}_{p_l}=&2u_l\delta_{ij}-\delta_{il}u_j-\delta_{jl}u_i.
\end{split}
\end{align}
By the definition of $v$, we have
\begin{align}\label{3vvk}
\begin{split}
v v_k=&u_1u_{1k}.
\end{split}
\end{align}
Hence, from\eqref{3Deltau}, we have
\begin{align}\label{3aijuijkb}
\begin{split}
\sum_{1\leq i,j\leq n}a^{ij}u_{ijk}
=&-2u_1u_{1k}\Delta u+2u_1\sum_{1\leq i\leq n}u_{1i}u_{ik}+v^3D_{k}f+3fvu_1u_{1k},\\
=&\frac{2u_1}{v^2}u_{11}u_{1k}+2u_1\sum_{2\leq i\leq n}u_{1i}u_{ik}+v^3D_{k}f+fvu_1u_{1k}.
\end{split}
\end{align}
By \eqref{3aijuijkb}, we get
\begin{align}\label{3I13}
\begin{split}
I_{13}
=&\frac{4 u_1^2}{v^2} u_{11}\sum_{1\leq k\leq n}c^{k1}u_{1k}
+4u_1^2\sum_{2\leq i\leq n}u_{1i}\sum_{1\leq k\leq n}c^{k1}u_{ki}+2fu_1^2v\sum_{1\leq k\leq n}c^{k1}u_{1k}\\&+2u_1v^3\sum_{1\leq k\leq n}c^{k1}D_{k}f.
\end{split}
\end{align}

For the term $I_{14}$:
\begin{align}\label{3I14}
\begin{split}
I_{14}
=&2u_{11}\sum_{1\leq k\leq n}c^{k1}u_{k1}+2v^2\sum_{2\leq i\leq n}u_{1i}\sum_{1\leq k\leq n}c^{k1}u_{ki}+2v^2\sum_{2\leq i\leq n}c^{1i}u_{1i}u_{ii}
\\
&+2u_{11}\sum_{2\leq i\leq n}c^{1i}u_{1i}+2\sum_{2\leq i,j\leq n}c^{ij}u_{1i}u_{1j}+2v^2\sum_{2\leq i\leq n}c^{ii}u_{ii}^2.
\end{split}
\end{align}

Combining  \eqref{3I11}, \eqref{3I12}, \eqref{3I13} and \eqref{3I14}, it follows that
\begin{align}\label{3I1c}
\begin{split}
I_1=&\frac{1}{|D'u|^2\log|D'u|^2}\bigg[(\frac{4 u_1^2}{v^2}+2)u_{11}\sum_{1\leq k\leq n}c^{k1}u_{k1}+(4u_1^2+2v^2)\sum_{2\leq i\leq n}u_{1i}\sum_{1\leq k\leq n}c^{k1}u_{ki}\\
&+2v^2\sum_{2\leq i\leq n}c^{1i}u_{1i}u_{ii}
+2u_{11}\sum_{2\leq i\leq n}c^{1i}u_{1i}+2\sum_{2\leq i,j\leq n}c^{ij}u_{1i}u_{1j}
+2fvu_1^2\sum_{1\leq k\leq n}c^{k1}u_{1k}\\
&+4(c^{11})_{1} u_1u_{11}+4 u_1\sum_{2\leq i\leq n}[(c^{1i})_{1}+v^2(c^{11})_{i}]u_{1i}+2v^2\sum_{2\leq i\leq n}c^{ii}u_{ii}^2+4u_1v^2\sum_{2\leq i\leq n}(c^{1i})_{i}u_{ii}\\
&+2u_1v^3\sum_{1\leq k\leq n}c^{k1}D_{k}f+\sum_{2\leq i\leq n}(c^{11})_{ii} u^4_1+\sum_{1\leq i\leq n}(c^{11})_{ii} u^2_1\bigg].
\end{split}
\end{align}

Inserting \eqref{3I2b} and \eqref{3I1c} into \eqref{3I1a}, we  can obtain the following formula
\begin{align}\label{3aijvarphiijb}
\begin{split}
0\geq \sum_{1\leq i,j\leq n}a^{ij}\varphi_{ij}=: Q_1+Q_2+Q_3,
\end{split}
\end{align}
where $Q_1$ contains all the quadratic terms of $u_{ij}$; $Q_2$  is the term which contains all linear terms of $u_{ij}$;  and   the remaining terms are denoted by $Q_3$.
Then we have
\begin{align}\label{3Q1a}
Q_1=&\frac{1}{|D'u|^2\log|D'u|^2}\big[ (\frac{4 u_1^2}{v^2}+2)u_{11}\sum_{1\leq k\leq n}c^{k1}u_{k1}+(4u_1^2+2v^2)\sum_{2\leq i\leq n}u_{1i}\sum_{1\leq k\leq n}c^{k1}u_{ki}\notag\\
&\qquad\qquad\qquad\qquad+2v^2\sum_{2\leq i\leq n}c^{1i}u_{1i}u_{ii}
+2u_{11}\sum_{2\leq i\leq n}c^{1i}u_{1i}\notag\\
&\qquad\qquad\qquad\qquad+2\sum_{2\leq i,j\leq n}c^{ij}u_{1i}u_{1j}
+2v^2\sum_{2\leq i\leq n}c^{ii}u_{ii}^2\big].
\end{align}
The linear terms of $u_{ij}$ are
\begin{align}\label{3Q2a}
Q_2=&\frac{1}{|D'u|^2\log|D'u|^2}\big[
2fvu_1^2\sum_{1\leq k\leq n}c^{k1}u_{1k}+4(c^{11})_{1} u_1u_{11}+4 u_1\sum_{2\leq i\leq n}(c^{1i})_{1}u_{1i}\notag\\
&\qquad\qquad\qquad\quad\,\,+4 u_1v^2\sum_{2\leq i\leq n}(c^{11})_{i}u_{1i}
+4u_1v^2\sum_{2\leq i\leq n}(c^{1i})_{i}u_{ii}\big],
\end{align}
and  the
remaining terms are
\begin{align}
Q_3=&I_2+\frac{1}{|D'u|^2\log|D'u|^2}\big[\sum_{2\leq i\leq n}(c^{11})_{ii} u^4_1+\sum_{1\leq i\leq n}(c^{11})_{ii} u^2_1+2u_1v^3\sum_{1\leq k\leq n}c^{k1}D_kf\big]\notag \\
=&I_2+\frac{1}{|D'u|^2\log|D'u|^2}\big[\sum_{2\leq i\leq n}(c^{11})_{ii} u^4_1+\sum_{1\leq i\leq n}(c^{11})_{ii} u^2_1+2c^{11}f_uu_1^2v^3\notag \\&
\hspace{120pt}+2u_1v^3\sum_{1\leq k\leq n}c^{k1}f_{x_k}\big].\label{3Q3a}
\end{align}
From the estimate on $I_2$ in \eqref{3I2c}, we have
\begin{align}\label{3Q3b}
Q_3\geq \sqrt{n}\alpha_0f v^3-2(n+c^{11}) \alpha_0^2u_1^2\log u_1-C_6u_1^2,
\end{align}
in the computation of $Q_3$, we use the relation $D_kf=f_u u_k+f_{x_k}$ and $f_u\geq 0$, where $C_6$ is a positive constant which depends only on $n, \Omega, M_0, \mu_0, L_1, L_2$.

{\bf Step 2:} In this step we shall treat the terms $Q_1, Q_2$ using the first order derivative condition
 $$\varphi_i(x_0)=0,$$
  and let
  \begin{align}\label{3A}
A=|D'u|^2\log|D'u|^2.
\end{align}

By \eqref{3varphii=0}  and  \eqref{3|D'u|^2i},  we have
\begin{align}\label{3ck1uki=1}
\begin{split}
\sum_{1\leq k\leq n}c^{k1}u_{ki}=&-\frac{h'}{2}\frac{u_i}{u_1}|D'u|^2\log|D'u|^2-\frac{g'\gamma^i}{2}\frac{|D'u|^2\log|D'u|^2}{u_1}-\frac{(c^{11})_i}{2} u_1\\
=&-\frac{h'}{2}\frac{u_i}{u_1}A-\frac{g'\gamma^i}{2}\frac{A}{u_1}-\frac{(c^{11})_i}{2} u_1,\qquad i=1,2,\ldots,n.
\end{split}
\end{align}
Using \eqref{3ck1uki=1},  we get
\begin{align}\label{3ck1uk1}
\begin{split}
\sum_{1\leq k\leq n} c^{k1}u_{k1}
=-\frac{h'}{2} A-\frac{g'\gamma^1}{2}\frac{A}{u_1}-\frac{(c^{11})_1}{2} u_1,
\end{split}
\end{align}
and
\begin{align}\label{3ck1uki=2}
\begin{split}
\sum_{1\leq k\leq n} c^{k1}u_{ki}
=-\frac{g'\gamma^i}{2}\frac{A}{u_1}-\frac{(c^{11})_i}{2} u_1,
 \qquad\qquad\quad i=2,3,\ldots,n.
\end{split}
\end{align}
Through \eqref{3ck1uki=2}  and the choice of the coordinate at $x_0$, we have
\begin{align}\label{3u1i}
\begin{split}
u_{1i}
=-\frac{c^{1i}}{c^{11}}u_{ii}-\frac{g'\gamma^i}{2c^{11}}\frac{A}{u_1}-\frac{(c^{11})_i}{2c^{11}} u_1,\quad\quad
\ \qquad i=2,3,\ldots,n.
\end{split}
\end{align}
Using \eqref{3ck1uk1}  and \eqref{3u1i},  it follows that
\begin{align}\label{3u11b}
\begin{split}
u_{11}
=&\sum_{2\leq i\leq n}\frac{(c^{1i})^2}{(c^{11})^2}u_{ii}-\frac{h'}{2c^{11}} A-\frac{g'\gamma^1}{c^{11}}\frac{A}{u_1}
+\frac{u_1}{2(c^{11})^2}\sum_{2\leq i\leq n}c^{1i}(c^{11})_i-\frac{(c^{11})_1}{2c^{11}}u_1\\
=&\sum_{2\leq i\leq n}\frac{(c^{1i})^2}{(c^{11})^2}u_{ii}-\frac{h'}{2c^{11}} A-\frac{g'\gamma^1}{c^{11}}\frac{A}{u_1}+bu_1,
\end{split}
\end{align}
 where we have let $$b=\frac{1}{2(c^{11})^2}\sum_{2\leq i\leq n}c^{1i}(c^{11})_i-\frac{(c^{11})_1}{2c^{11}}.$$
By \eqref{3u11a}  and \eqref{3u11b},  we have
\begin{align}\label{3uii=2}
\begin{split}
\sum_{2\leq i\leq n}\big[(c^{11})^2v^2+(c^{1i})^2\big]u_{ii}
=(c^{11})^2f v^3+\frac{c^{11}h'}{2} A+c^{11}g'\gamma^1\frac{A}{u_1}-(c^{11})^2bu_1.
\end{split}
\end{align}

Now we use the formulas \eqref{3ck1uk1}-\eqref{3u11b}  to treat  each term in $Q_1, Q_2$.  At first, we treat  the first five terms of $Q_1$ in \eqref{3Q1a}, and get \eqref{3Q1a1}-\eqref{3Q1a5}.

By \eqref{3ck1uk1} and  \eqref{3u11b}, we have
\begin{align}\label{3Q1a1}
\begin{split}
&(\frac{4 u_1^2}{v^2}+2)u_{11}\sum_{1\leq k\leq n}c^{k1}u_{k1}\\
=&(\frac{4 u_1^2}{v^2}+2)\big[\sum_{2\leq i\leq n}\frac{(c^{1i})^2}{(c^{11})^2}u_{ii}-\frac{h'}{2c^{11}} A-\frac{g'\gamma^1}{c^{11}}\frac{A}{u_1}+bu_1\big]
\big[-\frac{h'}{2} A-\frac{g'\gamma^1}{2}\frac{A}{u_1}-\frac{(c^{11})_1}{2} u_1\big]\\
=&-(\frac{2u_1^2}{v^2}+1)\big[h' A+g'\gamma^1\frac{A}{u_1}+(c^{11})_1 u_1\big]\sum_{2\leq i\leq n}\frac{(c^{1i})^2}{(c^{11})^2}u_{ii}
+(\frac{ u_1^2}{v^2}+\frac{1}{2})\frac{h'^2}{c^{11}} A^2\\&+(\frac{ u_1^2}{v^2}+\frac{1}{2})\frac{3h'g'\gamma^1}{c^{11}}\frac{A^2}{u_1}
+(\frac{ u_1^2}{v^2}+\frac{1}{2})[\frac{(c^{11})_1}{c^{11}}-2b]h'u_1A+(\frac{ 2u_1^2}{v^2}+1)\frac{g'^2(\gamma^1)^2}{c^{11}}\frac{A^2}{u_1^2}\\&
+(\frac{ 2u_1^2}{v^2}+1)g'\gamma^1[\frac{(c^{11})_1}{c^{11}}-b]A-(\frac{ 2u_1^2}{v^2}+1)(c^{11})_1bu_1^2.
\end{split}
\end{align}
 From \eqref{3ck1uki=2} and \eqref{3u1i}, we get
\begin{align}\label{3Q1a2}
\begin{split}
&(4u_1^2+2v^2)\sum_{2\leq i\leq n}u_{1i}\sum_{1\leq k\leq n}c^{k1}u_{ki}\\
=&(4u_1^2+2v^2)\sum_{2\leq i\leq n}\big[-\frac{c^{1i}}{c^{11}}u_{ii}-\frac{g'\gamma^i}{2c^{11}}\frac{A}{u_1}-\frac{(c^{11})_i}{2c^{11}} u_1\big]
\big[-\frac{g'\gamma^i}{2}\frac{A}{u_1}-\frac{(c^{11})_i}{2} u_1\big]\\
=&(2u_1^2+v^2)\frac{A}{u_1}\frac{g'}{c^{11}}\sum_{2\leq i\leq n}c^{1i}\gamma^iu_{ii}+(2u_1^2+v^2)u_1\sum_{2\leq i\leq n}\frac{c^{1i}}{c^{11}}(c^{11})_i u_{ii}+\frac{3g'^2}{2}A^2\\
&+\frac{g'}{c^{11}}\sum_{2\leq i\leq n}(c^{11})_i\gamma^i(2u_1^2+v^2)A+\frac{(2u_1^2+v^2)u_1^2}{2c^{11}}\sum_{2\leq i\leq n}((c^{11})_i)^2+\frac{g'^2}{2}\frac{A^2}{u_1^2}.
\end{split}
\end{align}
From \eqref{3u1i} , we have
\begin{align}\label{3Q1a3}
\begin{split}
&2v^2\sum_{2\leq i\leq n}c^{1i}u_{1i}u_{ii}\\
=&2v^2\sum_{2\leq i\leq n}c^{1i}u_{ii}\big[-\frac{c^{1i}}{c^{11}}u_{ii}-\frac{g'\gamma^i}{2c^{11}}\frac{A}{u_1}-\frac{(c^{11})_i}{2c^{11}} u_1\big]\\
=&-\frac{2v^2}{c^{11}}\sum_{2\leq i\leq n}(c^{1i}u_{ii})^2-\frac{g'}{c^{11}}\frac{v^2A}{u_1}\sum_{2\leq i\leq n}c^{1i}\gamma^iu_{ii}-\frac{u_1v^2}{c^{11}}\sum_{2\leq i\leq n}c^{1i}(c^{11})_iu_{ii}.
\end{split}
\end{align}
By \eqref{3u1i}  and \eqref{3u11b},  it follows that
\begin{align}\label{3Q1a4}
\begin{split}
&2u_{11}\sum_{2\leq i\leq n}c^{1i}u_{1i}\\=&2\big[\sum_{2\leq j\leq n}\frac{(c^{1j})^2}{(c^{11})^2}u_{jj}-\frac{h'}{2c^{11}} A-\frac{g'\gamma^1}{c^{11}}\frac{A}{u_1}+bu_1\big]\sum_{2\leq i\leq n}c^{1i}\big[-\frac{c^{1i}}{c^{11}}u_{ii}-\frac{g'\gamma^i}{2c^{11}}\frac{A}{u_1}-\frac{(c^{11})_i}{2c^{11}} u_1\big]\\
=&-\frac{2}{(c^{11})^3}[\sum_{2\leq i\leq n}(c^{1i})^2u_{ii}]^2+\big[h'A+3g'\gamma^1\frac{A}{u_1}-4c^{11}bu_1-(c^{11})_1u_1\big]\sum_{2\leq i\leq n}\frac{(c^{1i})^2}{(c^{11})^2}u_{ii}\\&-\frac{h'g'\gamma^1}{2c^{11}}\frac{A^2}{u_1}+\frac{h'}{2(c^{11})^2}\sum_{2\leq i\leq n}c^{1i}(c^{11})_i Au_1-\frac{g'^2(\gamma^1)^2}{c^{11}}\frac{A^2}{u_1^2}+g'\gamma^1\big[3b+\frac{(c^{11})_1}{c^{11}}\big]A\\&-b[2bc^{11}+(c^{11})_1]u^2_1.
\end{split}
\end{align}
Again by \eqref{3u1i}  and \eqref{2.2},  we get
\begin{align}\label{3Q1a5}
\begin{split}
&2\sum_{2\leq i,j\leq n}c^{ij}u_{1i}u_{1j}\\
=&2\sum_{2\leq i,j\leq n}c^{ij}\big[-\frac{c^{1i}}{c^{11}}u_{ii}-\frac{g'\gamma^i}{2c^{11}}\frac{A}{u_1}-\frac{(c^{11})_i}{2c^{11}} u_1\big]\big[-\frac{c^{1j}}{c^{11}}u_{jj}-\frac{g'\gamma^j}{2c^{11}}\frac{A}{u_1}-\frac{(c^{11})_j}{2c^{11}} u_1\big]\\
=&\frac{2}{(c^{11})^2}\sum_{2\leq i,j\leq n}c^{ij}c^{1i}c^{1j}u_{ii}u_{jj}-\big[\frac{2g'(\gamma^1)^3}{(c^{11})^2}\frac{A}{u_1}-\frac{2\gamma^1u_1}{(c^{11})^2}\sum_{2\leq j\leq n}\gamma^j(c^{11})_j\big]\sum_{2\leq i\leq n}(\gamma^i)^2u_{ii}\\&+\frac{2u_1}{(c^{11})^2}\sum_{2\leq i\leq n}c^{1i}(c^{11})_iu_{ii} +\frac{(1-c^{11})g'^2}{2c^{11}}\frac{A^2}{u_1^2}
+\frac{(1-c^{11})g'A}{(c^{11})^2 }\sum_{2\leq i\leq n}\gamma^i(c^{11})_i \\&+\frac{1}{2(c^{11})^2}\sum_{2\leq i,j\leq n}c^{ij}(c^{11})_i(c^{11})_j u_1^2.
\end{split}
\end{align}

Now we  treat  the first four terms of $Q_2$ in \eqref{3Q2a}, and get \eqref{3Q2a1}-\eqref{3Q2a4}.

From \eqref{3ck1uk1},  we have
\begin{align}\label{3Q2a1}
2fvu_1^2\sum_{1\leq k\leq n}c^{k1}u_{1k}
=&2fvu_1^2\big[-\frac{h'}{2} A-\frac{g'\gamma^1}{2}\frac{A}{u_1}-\frac{(c^{11})_1}{2} u_1\big]\notag\\
=&-h'fAvu_1^2-fg'\gamma^1Avu_1-(c^{11})_1fvu_1^3.
\end{align}
By \eqref{3u11b},   we obtain
\begin{align}\label{3Q2a2}
\begin{split}
4(c^{11})_{1} u_1u_{11}=&4(c^{11})_{1} u_1\big[\sum_{2\leq i\leq n}\frac{(c^{1i})^2}{(c^{11})^2}u_{ii}-\frac{h'}{2c^{11}} A-\frac{g'\gamma^1}{c^{11}}\frac{A}{u_1}+bu_1\big]\\
=&4(c^{11})_{1} u_1\sum_{2\leq i\leq n}\frac{(c^{1i})^2}{(c^{11})^2}u_{ii}-\frac{2(c^{11})_{1}}{c^{11}}h' Au_1-\frac{4g'\gamma^1}{c^{11}}(c^{11})_{1}A+4(c^{11})_{1} bu_1^2.
\end{split}
\end{align}
From \eqref{3u1i}, we have
\begin{align}\label{3Q2a3}
\begin{split}
&4 u_1\sum_{2\leq i\leq n}(c^{1i})_{1}u_{1i}\\
=&4 u_1\sum_{2\leq i\leq n}(c^{1i})_{1}\big[-\frac{c^{1i}}{c^{11}}u_{ii}-\frac{g'\gamma^i}{2c^{11}}\frac{A}{u_1}-\frac{(c^{11})_i}{2c^{11}} u_1\big]\\
=&-\frac{4u_1}{c^{11}}\sum_{2\leq i\leq n}(c^{1i})_{1}c^{1i}u_{ii}-\frac{2g'}{c^{11}}\sum_{2\leq i\leq n}(c^{1i})_{1}\gamma^iA-\frac{2}{c^{11}}\sum_{2\leq i\leq n}(c^{1i})_{1}(c^{11})_i u_1^2,
\end{split}
\end{align}
and
\begin{align}\label{3Q2a4}
\begin{split}
&4 u_1v^2\sum_{2\leq i\leq n}(c^{11})_{i}u_{1i}\\
=&4 u_1v^2\sum_{2\leq i\leq n}(c^{11})_{i}\big[-\frac{c^{1i}}{c^{11}}u_{ii}-\frac{g'\gamma^i}{2c^{11}}\frac{A}{u_1}-\frac{(c^{11})_i}{2c^{11}} u_1\big]\\
=&-\frac{4 u_1v^2}{c^{11}}\sum_{2\leq i\leq n}(c^{11})_{i}c^{1i}u_{ii}-\frac{2g'}{c^{11}}\sum_{2\leq i\leq n}(c^{11})_{i}\gamma^iAv^2-\frac{2}{c^{11}} \sum_{2\leq i\leq n}((c^{11})_{i})^2u_1^2v^2.
\end{split}
\end{align}
We treat the term $Q_1$ using the relations \eqref{3Q1a1}-\eqref{3Q1a5}, and  use the formulas \eqref{3Q2a1}-\eqref{3Q2a4}  to treat the term $Q_2$. By the formula on $Q_3$ in \eqref{3Q3a},  we can get the following new formula of \eqref{3aijvarphiijb},
\begin{align}\label{3aijvarphiijc}
0\geq \sum_{1\leq i,j\leq n}a^{ij}\varphi_{ij}=: J_1+J_2,
\end{align}
where $J_1$  only contains the terms with $u_{ii}$ , the other terms belong to $J_2$. We can write
\begin{align}\label{3J1a}
J_1=:\frac{1}{A}\big[J_{11}+J_{12}\big],
\end{align}
here $J_{11}$  contains the quadratic terms of $u_{ii}\,(i\geq2)$, and  $J_{12}$ is the term including linear terms of $u_{ii}\,(i\geq2)$. It follows that
\begin{align}\label{3J11a}
\begin{split}
J_{11}=&2v^2\sum_{2\leq i\leq n}c^{ii}u_{ii}^2-\frac{2v^2}{c^{11}}\sum_{2\leq i\leq n}(c^{1i}u_{ii})^2-\frac{2}{(c^{11})^3}[\sum_{2\leq i\leq n}(c^{1i})^2u_{ii}]^2\\&+\frac{2}{(c^{11})^2}\sum_{2\leq i,j\leq n}c^{ij}c^{1i}c^{1j}u_{ii}u_{jj}\\
=&\frac{2}{(c^{11})^3}\big[\sum_{2\leq i\leq n}d_ie_iu_{ii}^2+2\sum_{2\leq i<j\leq n}c^{ij}c^{1i}c^{1j}u_{ii}u_{jj}\big],
\end{split}
\end{align}
where
\begin{align}
d_i=&(c^{11})^2v^2+(c^{1i})^2=(c^{11})^2u_1^2+(c^{11})^2+(c^{1i})^2, \quad i=2,3,\ldots,n,\label{3di}\\
e_i=&c^{11}c^{ii}-(c^{1i})^2=1-(\gamma^1)^2-(\gamma^i)^2,\quad i=2,3,\ldots,n.\label{3ei}
\end{align}
And \begin{align}\label{3J12a}
\begin{split}
J_{12}=&\big[-\frac{2g'\gamma^1}{c^{11}}Au_1-\frac{2(\gamma^1)^2h'}{(c^{11})^2}\frac{Au_1^2}{v^2}
-\frac{2g'(\gamma^1)^3}{(c^{11})^2}\frac{A}{u_1}
+\frac{2\gamma^1}{(c^{11})^2}\sum_{2\leq j\leq n}c^{1j}(c^{11})_ju_1\\&\,\,\,-\frac{4b(\gamma^1)^2}{c^{11}}u_1+\frac{2g'(\gamma^1)^3}{(c^{11})^2}\frac{A}{v^2u_1}
+\frac{2(\gamma^1)^2(c^{11})_1}{(c^{11})^2}\frac{u_1}{v^2}\big]\sum_{2\leq i\leq n} (\gamma^i)^2u_{ii}\\&+4u_1v^2\sum_{2\leq i\leq n}(c^{1i})_{i}u_{ii}-\frac{4u_1}{c^{11}}\sum_{2\leq i\leq n}c^{1i}(c^{1i})_1u_{ii}\\&
-\big[\frac{2u_1^3}{c^{11}}+\frac{4c^{11}-2}{(c^{11})^2}u_1\big]\sum_{2\leq i\leq n}c^{1i}(c^{11})_i u_{ii}.
\end{split}
\end{align}

We write other terms as $J_2$, then
\begin{align*}
J_2=&Q_3-h'fvu_1^2+(\frac{ u_1^2}{v^2}+\frac{1}{2})\frac{h'^2}{c^{11}}A+\frac{3g'^2}{2}A
-fg'\gamma^1vu_1+\frac{g'}{c^{11}}\sum_{2\leq i\leq n}(c^{11})_i\gamma^i(u_1^2-1)\notag\\&
-\frac{1}{c^{11}} \sum_{2\leq i\leq n}((c^{11})_{i})^2\frac{(v^2+2)u_1^2}{A}-(c^{11})_1f\frac{v u_1^3}{A}
+(\frac{ 3u_1^2}{v^2}+1)\frac{h'g'\gamma^1}{c^{11}}\frac{A}{u_1}-\frac{2(c^{11})_{1}}{c^{11}}h' u_1\notag\\&
+(\frac{ u_1^2}{v^2}+\frac{1}{2})[\frac{(c^{11})_1}{c^{11}}-2b]h'u_1
+\frac{h'}{2(c^{11})^2}\sum_{2\leq i\leq n}c^{1i}(c^{11})_i u_1+\frac{g'^2}{2}\frac{A}{u_1^2}-\frac{g'^2(\gamma^1)^2}{c^{11}}\frac{A}{u_1^2}\notag\\&
+(\frac{2u_1^2}{v^2}+1)\frac{g'^2(\gamma^1)^2}{c^{11}}\frac{A}{u_1^2}
+\frac{(1-c^{11})g'^2}{2c^{11}}\frac{A}{u_1^2}+g'\gamma^1\big[3b+\frac{(c^{11})_1}{c^{11}}\big]-\frac{2g'}{c^{11}}\sum_{2\leq i\leq n}(c^{1i})_{1}\gamma^i
\end{align*}
\begin{align}
&+\frac{(1-c^{11})g'}{(c^{11})^2}\sum_{2\leq i\leq n}\gamma^i(c^{11})_i
+(\frac{ 2u_1^2}{v^2}+1)g'\gamma^1[\frac{(c^{11})_1}{c^{11}}-b]-(\frac{ 2u_1^2}{v^2}+1)(c^{11})_1b\frac{u_1^2}{A}\notag
\\&-b[2bc^{11}+(c^{11})_1]\frac{u_1^2}{A}
+\frac{1}{2(c^{11})^2}\sum_{2\leq i,j\leq n}c^{ij}(c^{11})_i(c^{11})_j \frac{u_1^2}{A}
-\frac{4g'\gamma^1}{c^{11}}(c^{11})_{1}\notag\\
&+4(c^{11})_{1} b\frac{u_1^2}{A}-\frac{2}{c^{11}}\sum_{2\leq i\leq n}(c^{1i})_{1}(c^{11})_i \frac{u_1^2}{A}.\label{3J2a}
\end{align}

Using the formula on $Q_3$ in  \eqref{3Q3b} and $I_2$ in  \eqref{3I2c}, we get the following estimate on $J_2$,
\begin{align}
J_2\geq&-2(h'^2+c^{11}g'^2) u_1^2\log u_1+h'fv+\frac{3}{2}\frac{h'^2}{c^{11}}A+\frac{3g'^2}{2}A-C_7u_1^2\notag\\
\ge & [h'^2+c^{11}g'^2]u_1^2\log u_1-C_8u_1^2.\label{3J2b}
\end{align}

So if we use $h(u), g(d)$ in \eqref{3hg1}, then we have
 \begin{align}\label{3J2c}
\begin{split}
J_2 \geq (n+c^{11}) \alpha_0^2u_1^2\log u_1-C_8u_1^2,
\end{split}
\end{align}
where  $ C_7, C_8$  and the following $C_9, ..., C_{15}$ are  positive constants which only depend on $n, \Omega, \mu_0, M_0, L_1, L_2$.

{\bf Step 3:}
In this step, we concentrate on $J_1$. We first treat the terms $J_{11}$ and $J_{12}$  and obtain the formula \eqref{3J1b}, then we complete the proof of  Theorem~\ref{Thm1.1} through  Lemma~\ref{Lem4.5}.\\
 By \eqref{3uii=2}, we have
\begin{align}\label{3u22}
\begin{split}
u_{22}=&-\frac{1}{d_2}\sum_{3\leq i\leq n}d_iu_{ii}
+\frac{1}{d_2}\big[(c^{11})^2f v^3+\frac{c^{11}h'}{2} A+c^{11}g'\gamma^1\frac{A}{u_1}
-(c^{11})^2bu_1\big]\\
=&:-\frac{1}{d_2}\sum_{3\leq i\leq n}d_iu_{ii}
+\frac{D}{d_2},
\end{split}
\end{align}
where we have let
\begin{align}\label{3D}
D=(c^{11})^2f v^3+\frac{c^{11}h'}{2} A+c^{11}g'\gamma^1\frac{A}{u_1}
-(c^{11})^2bu_1.
\end{align}

We first treat the term $J_{11}$:
using \eqref{3D} to simplify \eqref{3J11a}, we get
\begin{align}\label{3J11b}
\begin{split}
J_{11}=&\frac{2}{(c^{11})^3d_2}\big[\sum_{3\leq i\leq n}b_{ii}u_{ii}^2+2\sum_{3\leq i< j\leq n}b_{ij}u_{ii}u_{jj}
-2e_2D\sum_{3\leq i\leq n}d_{i}u_{ii}\\&\qquad\qquad\,\,-2(\gamma^2)^2D\sum_{3\leq i\leq n}(c^{1i})^2u_{ii}+e_2D^2\big],
\end{split}
\end{align}
where
\begin{align}\label{3bij}
\begin{split}
b_{ii}=&e_2d_i^2+e_id_id_2-2c^{12}c^{2i}c^{1i}d_i
=(c^{11})^4(e_2+e_i)v^4+A_{1i}v^2+A_{2i},\quad i\geq 3\\
b_{ij}=&e_2d_id_j+d_2c^{ij}c^{1i}c^{1j}-c^{12}c^{1i}c^{2i}d_j-c^{12}c^{1j}c^{2j}d_i\\
=&(c^{11})^4e_2v^4+G_{ij}v^2+\hat{G}_{ij},
\qquad\qquad\qquad\qquad i\neq j,\,i,j\geq 3,
\end{split}
\end{align}
and
\begin{align}\label{3Aij}
\begin{split}
A_{1i}=&(c^{11})^2\big[(c^{1i})^2(e_2+e_i)+c^{11}\big((c^{1i})^2+(c^{12})^2\big)\big], \\
A_{2i}=&c^{11}(c^{1i})^2\big[(c^{1i})^2+(c^{12})^2\big], \\
G_{ij}=&c^{11}\big((c^{1i})^2+(c^{1j})^2\big)+c^{ij}c^{1i}c^{1j}, \\
\hat{G}_{ij}=&c^{11}(c^{1i})^2(c^{1j})^2.
\end{split}
\end{align}
Now we simplify the terms in $J_{12}$: by  \eqref{3u22},  we can rewrite \eqref{3J12a} as
\begin{align}\label{3J12b}
\begin{split}
J_{12}=&\big[-\frac{2g'\gamma^1}{c^{11}}Au_1-\frac{2(\gamma^1)^2h'}{(c^{11})^2}\frac{Au_1^2}{v^2}
-\frac{2g'(\gamma^1)^3}{(c^{11})^2}\frac{A}{u_1}
+\frac{2\gamma^1}{(c^{11})^2}\sum_{2\leq j\leq n}c^{1j}(c^{11})_ju_1\\&\,\,\,-\frac{4b(\gamma^1)^2}{c^{11}}u_1+\frac{2g'(\gamma^1)^3}{(c^{11})^2}\frac{A}{v^2u_1}
+\frac{2(\gamma^1)^2(c^{11})_1}{(c^{11})^2}\frac{u_1}{v^2}\big]\sum_{3\leq i\leq n}[(\gamma^i)^2-\frac{d_i}{d_2}(\gamma^2)^2]u_{ii}\\&+4u_1v^2\sum_{3\leq i\leq n}[(c^{1i})_{i}-\frac{d_i}{d_2}(c^{12})_{2}]u_{ii}-\frac{4u_1}{c^{11}}\sum_{3\leq i\leq n}[c^{1i}(c^{1i})_1-\frac{d_i}{d_2}c^{12}(c^{12})_1]u_{ii}\\&
-\big[\frac{2u_1^3}{c^{11}}+\frac{4c^{11}-2}{(c^{11})^2}u_1\big]\sum_{3\leq i\leq n}[c^{1i}(c^{11})_i-\frac{d_i}{d_2}c^{12}(c^{11})_2]u_{ii}\\
&+\big[-\frac{2g'\gamma^1}{c^{11}}Au_1-\frac{2(\gamma^1)^2h'}{(c^{11})^2}\frac{Au_1^2}{v^2}
-\frac{2g'(\gamma^1)^3}{(c^{11})^2}\frac{A}{u_1}
+\frac{2\gamma^1}{(c^{11})^2}\sum_{2\leq j\leq n}c^{1j}(c^{11})_ju_1\\&\qquad-\frac{4b(\gamma^1)^2}{c^{11}}u_1+\frac{2g'(\gamma^1)^3}{(c^{11})^2}\frac{A}{v^2u_1}
+\frac{2(\gamma^1)^2(c^{11})_1}{(c^{11})^2}\frac{u_1}{v^2}\big]\frac{(\gamma^2)^2D}{d_2}\\&
+4(c^{12})_{2}\frac{u_1v^2D}{d_2}-4c^{12}(c^{12})_1\frac{u_1D}{c^{11}d_2}-\big[\frac{2u_1^3}{c^{11}}
+c^{12}(c^{11})_2\frac{4c^{11}-2}{(c^{11})^2}u_1\big]\frac{D}{d_2}.
\end{split}
\end{align}
Using \eqref{3J11b} and \eqref{3J12b} to treat \eqref{3J1a},  we have
\begin{align}\label{3J1b}
\begin{split}
J_{1}=&\frac{2}{Ad_2(c^{11})^3}\big[\sum_{3\leq i\leq n}b_{ii}u_{ii}^2+2\sum_{3\leq i< j\leq n}b_{ij}u_{ii}u_{jj}-u_1^5\log u_1\sum_{3\leq i\leq n}b_iu_{ii}\\&\qquad\qquad\quad+\sum_{3\leq i\leq n}K_iu_{ii}\big]+ R,
\end{split}
\end{align}
where
\begin{align}\label{3bi}
b_i=&2(c^{11})^5g'\gamma^1(e_2-e_i),
\end{align}
 and
\begin{align}\label{3Ki}
K_i=&-2e_2Dd_{i}-2(\gamma^2)^2D(c^{1i})^2-(c^{11})^3(\log c^{11})g'\gamma^1u_1^5(e_2-e_i)-(c^{11})^4g'\gamma^1Au_1(e_2-e_i)\notag\\&
+\frac{(c^{11})^3}{2}d_2\big[-\frac{2(\gamma^1)^2h'}{(c^{11})^2}\frac{Au_1^2}{v^2}
-\frac{2g'(\gamma^1)^3}{(c^{11})^2}\frac{A}{u_1}
+\frac{2\gamma^1}{(c^{11})^2}\sum_{2\leq j\leq n}c^{1j}(c^{11})_ju_1\notag
\\&-\frac{4b(\gamma^1)^2}{c^{11}}u_1+\frac{2g'(\gamma^1)^3}{(c^{11})^2}\frac{A}{v^2u_1}
+\frac{2(\gamma^1)^2(c^{11})_1}{(c^{11})^2}\frac{u_1}{v^2}\big][(\gamma^i)^2-\frac{d_i}{d_2}(\gamma^2)^2]\notag\\&+4u_1v^2\frac{(c^{11})^3}{2}d_2\sum_{3\leq i\leq n}[(c^{1i})_{i}-\frac{d_i}{d_2}(c^{12})_{2}]-\frac{4u_1}{c^{11}}\frac{(c^{11})^3}{2}d_2 [c^{1i}(c^{1i})_1-\frac{d_i}{d_2}c^{12}(c^{12})_1]\notag\\&
-\frac{(c^{11})^3}{2}d_2\big[\frac{2u_1^3}{c^{11}}+\frac{4c^{11}-2}{(c^{11})^2}u_1\big][c^{1i}(c^{11})_i-\frac{d_i}{d_2}c^{12}(c^{11})_2]
\end{align}
we also have let
\begin{align*}
\begin{split}
R=&\frac{2e_2D^2}{(c^{11})^3d_2A}+\big[-\frac{2g'\gamma^1}{c^{11}}u_1-\frac{2(\gamma^1)^2h'}{(c^{11})^2}\frac{u_1^2}{v^2}
-\frac{2g'(\gamma^1)^3}{(c^{11})^2}\frac{1}{u_1}
+\frac{2\gamma^1}{(c^{11})^2}\sum_{2\leq j\leq n}c^{1j}(c^{11})_j\frac{u_1}{A}\\&\,\,\,-\frac{4b(\gamma^1)^2}{c^{11}}\frac{u_1}{A}+\frac{2g'(\gamma^1)^3}{(c^{11})^2}\frac{1}{v^2u_1}
+\frac{2(\gamma^1)^2(c^{11})_1}{(c^{11})^2}\frac{u_1}{v^2A}\big]\frac{(\gamma^2)^2D}{d_2}\\&
+4(c^{12})_{2}\frac{u_1v^2D}{Ad_2}-4c^{12}(c^{12})_1\frac{u_1D}{c^{11}Ad_2}-\big[\frac{2u_1^3}{c^{11}}
+c^{12}(c^{11})_2\frac{4c^{11}-2}{(c^{11})^2}u_1\big]\frac{D}{Ad_2}.
\end{split}
\end{align*}

For $K_i$ and $R$, using the formulas on $D$ in \eqref{3D}; the formula of $A$ in \eqref{3A}; $ e_i; d_i$  in \eqref{3di}-\eqref{3ei}, and  $h(u), g(d)$ in \eqref{3hg1}, we have the following estimates
\begin{align}
K_i\leq& C_9u_1^5,\label{3Kia}\\
R\leq &C_{10}u_1^2.\label{3R}
\end{align}

Now we use  Lemma~\ref{Lem4.5}, if there is a sufficiently large positive constant $C_{11}$ such that
\begin{align}
|Du|(x_0) \ge C_{11},\label{3C11}
\end{align}
 then we have
\begin{align}\label{3J1c}
J_{1}
\geq &\frac{2}{Ad_2(c^{11})^3}{[-(n-2)(c^{11})^7g'^2(\gamma^1)^2u_1^6\log^2 u_1-C_{12} u_1^6\log u_1]}-C_{10}u_1^2,\notag\\
\ge & -(n-2)c^{11}(1-c^{11})g'^2u_1^2\log u_1-C_{13}u_1^2,
\end{align}
where we use the formulas  $(\gamma^1)^2= 1-c^{11}$, $d_2$ in \eqref{3di} and $A$ in \eqref{3A}.

Using the estimates on $J_1$ in \eqref{3J1c} and  $J_2$ in \eqref{3J2b}, from  \eqref{3aijvarphiijc} we obtain
\begin{align}\label{3aijvarphiijd}
0\geq &\sum_{1\leq i,j\leq n}a^{ij}\varphi_{ij}\notag\\
\geq & \big\{h'^2+[(c^{11})^2(n-2)-c^{11}(n-3)]g'^2\big\}u_1^2\log u_1-C_{14}u_1^2.
\end{align}
By the choice of $h(u), g(d)$ in \eqref{3hg1},
it follows that
\begin{align}\label{3aijvarphiije}
\begin{split}
0\geq &\sum_{1\leq i,j\leq n}a^{ij}\varphi_{ij}\\
\ge & \big\{n+[(c^{11})^2(n-2)-c^{11}(n-3)]\alpha_0^2\big\}u_1^2\log u_1-C_{14}u_1^2\\
\geq & 3u_1^2\log u_1-C_{14}u_1^2.
\end{split}
\end{align}

By \eqref{3c^{11}}, \eqref{3C11} and \eqref{3aijvarphiije}, there exists a positive constant
$C_{15}$ such that
 \begin{align}\label{3C15}
 |D'u|(x_0)\leq C_{15}.
 \end{align}

 So from  Case I, Case II,  and \eqref{3C15}, we have
 $$|D'u|(x_0)\leq C_{16}, \quad \quad x_0\in\Omega_{\mu_0}\bigcup\partial\Omega.$$

Since $\varphi(x)\leq\varphi(x_0),\quad \text{for} \,\forall x\in \Omega_{\mu_0}$, there exists $M_2$ such that
\begin{align}\label{3M2}
|Du|(x)\leq M_2, \quad in\quad\Omega_{\mu_0}\bigcup\partial\Omega,
\end{align}
where $M_2$ depends only on $n, \Omega, \mu_0, M_0,  L_1, L_2$.

So at last we get the following estimate
$$\sup_{\overline\Omega_{\mu_0}}|Du|\leq \max\{M_1, M_2\},$$
where the positive constant  $ M_1$ depends only on $n, \mu_0, M_0, L_1$; and $ M_2$ depends only on $n, \Omega, \mu_0, M_0, L_1, L_2$. Now we  complete the proof of Theorem~\ref{Thm1.1}.\qed

\section{ Some  Lemmas }
In this section, we prove the main Lemma~\ref{Lem4.5} and get the main estimate \eqref{3J1c}, which was used in last section to estimate  $J_1$ defined in  \eqref{3J1b}.

We first state a simple lemma on elementary symmetric function.
\begin{Lem}\label{Lem4.1}
Assume $e=(e_2, e_3,\ldots, e_n)$,  then for $i\geq 3$, we have
\begin{align}\label{lem4.1}
\sigma_{n-3}(e|i)(e_2-e_i)-\sum_{k\neq i, k\geq 3}\sigma_{n-3}(e|ik)(e_2-e_k)
=(n-1)\sigma_{n-2}(e|i)-\sigma_{n-2}(e).
\end{align}
\end{Lem}
{\bf Proof:}  When ~$i\geq 3$, we have,
\begin{align*}
\begin{split}
&\sigma_{n-3}(e|i)(e_2-e_i)-\sum_{k\neq i, k\geq 3}\sigma_{n-3}(e|ik)(e_2-e_k)\\
=&e_2\sigma_{n-3}(e|i)- e_i\sigma_{n-3}(e|i)-e_2\sum_{i\neq k, k\geq 3}\sigma_{n-3}(e|ik)+\sum_{k\neq i, k\geq 3}e_k\sigma_{n-3}(e|ik)\\
=&(n-2)\sigma_{n-2}(e|i)-e_i\sigma_{n-3}(e|i)\\
=&(n-1)\sigma_{n-2}(e|i)-\sigma_{n-2}(e).
\end{split}
\end{align*}
\qed

\begin{Lem}\label{Lem4.2} Let $a_i=(\gamma^i)^2, \gamma = (\gamma^1, \gamma^2,...,\gamma^n) \quad \text{is  a unit vector in} \quad R^n, \quad a=(a_2, a_3, \ldots, a_n)$, and $e=(e_2, e_3, \ldots, e_n), \quad e_i=\sigma_1(a|i)$, $i\geq 2$. Then the matrix
$E=(E_{ij})_{3\leq i,j\leq n}$ is positive definite, where $ E_{ij}=e_2+e_i\delta_{ij}$.
\end{Lem}
{\bf Proof:} We only need to prove that the following determination is positive.
\begin{align}\label{4.1}
\begin{split}
\det E=&\sigma_{n-2}(e)
=\sigma_{n-2}\big(\sigma_1(a|2), \sigma_1(a|3), \ldots, \sigma_1(a|n)\big)\\
=&\sum_{2\leq i_1<i_2<\cdots< i_{n-2}\leq n-2}\big(\sigma_1(a)-a_{i_1}\big)\big(\sigma_1(a)-a_{i_2}\big)\cdots \big(\sigma_1(a)-a_{i_{n-2}}\big)\\
=&\sum_{0\leq k\leq n-2}(-1)^k(n-k-1)[\sigma_1(a)]^{n-2-k}\sigma_k(a)\\
=&[\sigma_1(a)]^{n-2}+\sum_{2\leq k\leq n-2}(-1)^k(n-k-1)[\sigma_1(a)]^{n-2-k}\sigma_k(a),
\end{split}
\end{align}
Now we divide the following two cases, using the  Newton-MacLaurin inequality, then we get our conclusion.

Case 1:  if $ n=odd$
\begin{align}\label{4.2}
\begin{split}
&\sum_{2\leq k\leq n-2}(-1)^k(n-1-k)[\sigma_1(a)]^{n-2-k}\sigma_k(a)\\
=&\sum_{2\leq k\leq n-3,k=even}[k(\sigma_1(a))^{k-1}\sigma_{n-1-k}(a)-(k-1)(\sigma_1(a))^{k-2}\sigma_{n-k}(a)]\\
=&\sum_{2\leq k\leq n-3,k=even}[\sigma_1(a)]^{k-2}[k \sigma_1(a)\sigma_{n-1-k}(a)-(k-1)\sigma_{n-k}(a)]\\
\geq&\sum_{2\leq k\leq n-3,k=even}[\sigma_1(a)]^{k-2}[(n-1)(n-k)-(k-1)]\sigma_{n-k}(a)\\
\geq& 0.
\end{split}
\end{align}
Case 2: if $ n=even$
\begin{align}\label{4.3}
\begin{split}
&\sum_{2\leq k\leq n-2}(-1)^k(n-1-k)[\sigma_1(a)]^{n-2-k}\sigma_k(a)\\
=&\sum_{3\leq k\leq n-3,k=odd}[k(\sigma_1(a))^{k-1}\sigma_{n-1-k}(a)-(k-1)(\sigma_1(a))^{k-2}\sigma_{n-k}(a)]+\sigma_{n-2}(a)\\
\geq&\sum_{3\leq k\leq n-3,k=odd}[\sigma_1(a)]^{k-2}[k\sigma_1(a)\sigma_{n-1-k}(a)-(k-1)\sigma_{n-k}(a)]+\sigma_{n-2}(a)\\
\geq&\sum_{3\leq k\leq n-3,k=odd}[\sigma_1(a)]^{k-2}[(n-1)(n-k)-(k-1)]\sigma_{n-k}(a)+\sigma_{n-2}(a)\\
\geq& 0.
\end{split}
\end{align}
Since  $\sigma_1(a)=\sum_{2\leq i\leq n}a_i=c^{11}>0,$
it follows that
\begin{align}\label{4.4}
\det E=\sigma_{n-2}(e)\geq[\sigma_1(a)]^{n-2}>0.
\end{align}
then the matrix $E$ is positive definite.\qed

Now we prove  the main lemma.
\begin{Lem}\label{Lem4.5}
We define $(b_{ij})$  as in \eqref{3bij}, $d_i,e_i$ defined as in \eqref{3di}-\eqref{3ei}, $A_{1i}, A_{2i}, G_{ij}, \hat{G}_{ij}$ defined as in \eqref{3Aij}.
And we define $b_i$  as in \eqref{3bi}, $v^2 = 1 + u_1^2$ and $c^{11}\ge \frac{1}{C_4}.$  We study the following quadratic form
\begin{align}\label{4.6}
\begin{split}
Q(x_3,x_4,\ldots,x_n)=&\sum_{3\leq i\leq n}b_{ii}x_i^2
+2\sum_{3\leq i< j\leq n}b_{ij}x_i x_j
-u_1^5\log u_1\sum_{3\leq i\leq n}b_{i} x_i\\
&+\sum_{3\leq i\leq n}K_i x_i,
\end{split}
\end{align}
where $K_i$ defined in  \eqref{3Ki} and we have the estimate  \eqref{3Kia} for $K_i$.
Then there exists a sufficiently large positive constant $C_{16}$ which depends only on  $n, \Omega, \mu_0, M_0, L_1, L_2,$ such that if
\begin{align}\label{4.6a}
|Du|(x_0)=u_1(x_0) \ge C_{16},
\end{align}
then  the followings hold.

(I): The matrix $(b_{ij})$ is positive definite if and only only if the matrix\\
$(b^1_{ij})=  (E_{ij})=[e_2+e_i\delta_{ij}]$ is positive definite.\par

(II): We have
\begin{align}\label{4.7}
\begin{split}
Q(x_3,x_4,\ldots,x_n)
\geq &-(n-2)(c^{11})^7g'^2(\gamma^1)^2u_1^6\log^2 u_1-C_{17} u_1^6\log u_1,
\end{split}
\end{align}
where positive constant $C_{17}$ also depends only on  $n, \Omega, \mu_0, M_0, L_1, L_2.$

\end{Lem}
{\bf Proof:} Let $$ B=(b_{ij})=B_1+B_2, B_1=((c^{11})^4u_1^4b^1_{ij}), B_2=(O(u^2_1)\delta_{ij}),$$
We first  prove  (I):
\begin{align}\label{4.8}
\begin{split}
\sigma_k(B)=&\sigma_k(B_1+B_2)\\
=&\sigma_k(B_1)+\sigma_k(B_1, B_1, \ldots, B_1, B_2)\\&+\cdots+\sigma_k(B_1, B_2,\ldots, B_2, B_2)+\sigma_k(B_2)\\
=&(c^{11})^{4k}u_1^{4k}\sigma_k(b^1_{ij})+O(u_1^{4k-2})
\end{split}
\end{align}
so if $u_1$ is sufficiently large, then $\sigma_k(B)>0\Longleftrightarrow\sigma_k(b^1_{ij})>0$.

Now we prove   (II):
If  $B_1=((c^{11})^4u_1^4b^1_{ij})_{3\leq i,j\leq n}$  is positive definite, from the argument in (I), we get
\begin{align}\label{4.9}
B^{-1}=(B_1+B_2)^{-1}=B_1^{-1}(I+B_1^{-1}B_2)^{-1}=\frac{1}{(c^{11})^{4}u_1^{4}}(b^1_{ij})^{-1}\big(1+o(1)\big).
\end{align}
Then we have
\begin{align}\label{4.10}
\begin{split}
(b^1_{ij})^{-1}=&\left(\begin{array}{cccc}e_2+e_3&e_2&\cdots &e_2\\
e_2&e_2+e_4&\cdots &e_2\\
\vdots & \vdots & \vdots &\vdots \\
e_2&e_2&\cdots& e_2+e_n
\end{array}\right)^{-1}\\
=&\frac{1}{\sigma_{n-2}(e)}\left(\begin{array}{cccc}
\sigma_{n-3}(e|3)&-\sigma_{n-3}(e|34)&\cdots & -\sigma_{n-3}(e|3n)\\
-\sigma_{n-3}(e|43)&\sigma_{n-3}(e|4)&\cdots & -\sigma_{n-3}(e|4n)\\
\vdots & \vdots & \vdots &\vdots \\
-\sigma_{n-3}(e|n3)&-\sigma_{n-3}(e|n4)&\cdots & \sigma_{n-3}(e|n)
\end{array}\right)\\
=&:\frac{1}{\sigma_{n-2}(e)}\tilde{B}
\end{split}
\end{align}
where  $e=(e_2,e_3,\ldots,e_n)$.\\
Now we solve the following linear  algebra equation
\begin{align}\label{4.11}
\frac{\partial Q}{\partial x_k}=0,\quad k=3,4,\ldots,n.
\end{align}
 We assume $(\bar{x}_{3},\bar{x}_{4},\ldots,\bar{x}_{n})$ is the extreme point of the quadratic form $Q(x_3,x_4,\ldots,x_n)$.
From the definition of $b_{ij}, b_{i}, K_i$ in \eqref{3bij}, \eqref{3bi}, \eqref{3Ki}  and the estimate for $K_i$ in \eqref{3Kia}, using the formulas \eqref{4.9} and \eqref{4.10}, it follows that
\begin{align}\label{4.12}
\begin{split}
\left(\begin{array}{c}
\bar{x}_{3}\\
\bar{x}_{4}\\
\vdots  \\
\bar{x}_{n}
\end{array}\right)=&\frac{1}{2}u_1^5\log u_1B^{-1}
\left(\begin{array}{c}
b_{3}\\
b_{4}\\
\vdots\\
b_{n}
\end{array}\right)+ O(u_1^5)B^{-1}
\left(\begin{array}{c}
1\\
1\\
\vdots\\
1
\end{array}\right)\\
=&\frac{1}{2}u_1^5\log u_1B_1^{-1}\left(\begin{array}{c}
b_{3}\\
b_{4}\\
\vdots\\
b_{n}
\end{array}\right)+ O(u_1)
\left(\begin{array}{c}
1\\
1\\
\vdots\\
1
\end{array}\right)\\
=&\frac{c^{11}g'\gamma^1u_1\log u_1}{\sigma_{n-2}(e)}\tilde{B}
\left(\begin{array}{c}
e_2-e_3\\
e_2-e_4\\
\vdots\\
e_2-e_n
\end{array}\right)+ O(u_1)
\left(\begin{array}{c}
1\\
1\\
\vdots\\
1
\end{array}\right).\\
\end{split}
\end{align}
From Lemma \ref{Lem4.1}, we have for $i=3,4,\ldots,n,$
\begin{align}\label{4.13}
\begin{split}
\bar{x}_{i}=&\frac{c^{11}g'\gamma^1u_1\log u_1}{\sigma_{n-2}(e)}\big[\sigma_{n-3}(e|i)(e_2-e_i)-\sum_{k\neq i, k\geq 3}\sigma_{n-3}(e|ik)(e_2-e_k)\big]+ O(u_1)\\
=&\frac{c^{11}g'\gamma^1u_1\log u_1}{\sigma_{n-2}(e)}\big[(n-1)\sigma_{n-2}(e|i)-\sigma_{n-2}(e)\big]+ O(u_1).
\end{split}
\end{align}
It follows that we have the following minimum of the quadratic $Q$,
\begin{align}\label{4.14}
\begin{split}
&Q(\bar{x}_3, \bar{x}_4,\ldots, \bar{x}_n)\\
=&\frac{(c^{11})^6g'^2(\gamma^1)^2u_1^6\log^2 u_1}{\sigma^2_{n-2}(e)}\bigg\{\sum_{3\leq i\leq n}(e_2+e_i)\big[(n-1)\sigma_{n-2}(e|i)-\sigma_{n-2}(e)\big]^2
\\&+2e_2\sum_{3\leq i< j\leq n}\big[(n-1)\sigma_{n-2}(e|i)-\sigma_{n-2}(e)\big] \big[(n-1)\sigma_{n-2}(e|j)-\sigma_{n-2}(e)\big]\\&
-2\sigma_{n-2}(e)\sum_{3\leq i\leq n}(e_2-e_i)\big[(n-1)\sigma_{n-2}(e|i)-\sigma_{n-2}(e)\big]\bigg\}
+O(u_1^6\log u_1).
\end{split}
\end{align}
By the elementary computation, we have
\begin{align}\label{4.15}
\begin{split}
&\sum_{3\leq i\leq n}(e_2+e_i)\big[(n-1)\sigma_{n-2}(e|i)-\sigma_{n-2}(e)\big]^2
\\&+2e_2\sum_{3\leq i< j\leq n}\big[(n-1)\sigma_{n-2}(e|i)-\sigma_{n-2}(e)\big] \big[(n-1)\sigma_{n-2}(e|j)-\sigma_{n-2}(e)\big]\\&
-2\sigma_{n-2}(e)\sum_{3\leq i\leq n}(e_2-e_i)\big[(n-1)\sigma_{n-2}(e|i)-\sigma_{n-2}(e)\big]\\
=&e_2\big\{\sum_{3\leq i\leq n}\big[(n-1)\sigma_{n-2}(e|i)-\sigma_{n-2}(e)\big]\big\}^2\\&
-2e_2\sigma_{n-2}(e)\sum_{3\leq i\leq n}\big[(n-1)\sigma_{n-2}(e|i)-\sigma_{n-2}(e)\big]\\&
+\sum_{3\leq i\leq n}e_i\big[(n-1)\sigma_{n-2}(e|i)-\sigma_{n-2}(e)\big]^2\\&
+2\sigma_{n-2}(e)\sum_{3\leq i\leq n}e_i\big[(n-1)\sigma_{n-2}(e|i)-\sigma_{n-2}(e)\big]\\
=&-e_2\sigma_{n-2}^2(e)+(n-1)^2\sum_{3\leq i\leq n}e_i\sigma^2_{n-2}(e|i)-\sigma_{n-2}^2(e)\sigma_{1}(e|2)\\
=&[(n-1)^2\sigma_{n-1}(e)-\sigma_{1}(e)\sigma_{n-2}(e)]\sigma_{n-2}(e)\\
\geq &-\sigma_{1}(e)\sigma^2_{n-2}(e)\\
=&-(n-2)c^{11}\sigma^2_{n-2}(e).
\end{split}
\end{align}
Using \eqref{4.14} and \eqref{4.15},  we at last  get the following estimate
\begin{align}\label{4.16}
\begin{split}
Q(x_3, x_4,\ldots, x_n)\geq&Q(\bar{x}_3, \bar{x}_4,\ldots, \bar{x}_n)\\
\geq &-(n-2)(c^{11})^7g'^2(\gamma^1)^2u_1^6\log^2 u_1+O(u_1^6\log u_1).
\end{split}
\end{align}
 In this computation, the bounds in the coefficient on $O(u_1^6\log u_1), O(u_1^5), O(u_1)$ depend only on $n, \Omega,  M_0, \mu_0,  L_1, L_2$. Thus we complete this proof.\qed

\section{Proof of  Theorem~\ref{Thm1.2}.}
{\em Proof of  Theorem~\ref{Thm1.2}.}

As in the proof of  Theorem~\ref{Thm1.1}., let
$$\tilde{P}(x)=\log|D'u|^2e^{\sqrt{n}\tilde{\alpha}_0(M_0+1+u)}e^{\tilde{\alpha}_0d},$$ where we  have let
$$\tilde{\alpha}_0=4L_2\frac{\sqrt{1+a_0}}{\sqrt{1-b_0^2}}+2C_0+2, $$ which is a constant, and $$a_0=\max_{x\in\partial\Omega}\frac{2\psi^2}{1-\psi^2},$$
$C_0$ is also a positive constant depending only on $n,\Omega$.

Similarly,  let
$$\varphi(x)=\log \tilde{P}(x)=\log\log|D'u|^2+h(u)+g(d),$$
where in the capillary boundary value case, we choose
\begin{align}\label{5hg1}
h(u)=\sqrt{n}\tilde{\alpha}_0(M_0+1+u), \quad g(d)=\tilde{\alpha}_0 d.
\end{align}

 We assume that
$\varphi(x)$ attains its maximum at $x_0 \in \overline \Omega_{\mu_{0}}$, where $0<\mu_0<\mu_1$ is a sufficiently small number which we shall decide it  later.

{\bf Case 1.} If $x_0\in \partial \Omega$, we shall get the bound of $|D'u|(x_0)$.

Similar calculations to case I in the proof of Theorem~\ref{Thm1.1}., let $q=0$ in \eqref{3.11}, we get
\begin{align}\label{5varphigamma1}
\begin{split}
|D'u|^2\log |D'u|^2 \frac{\partial\varphi}{\partial\gamma}(x_0)
=&g'|D'u|^2\log |D'u|^2-2\sum_{1\leq i,k,l\leq n}c^{kl}u_iu_l(\gamma^i)_k\\&-\frac{2v}{1-\psi^2}\sum_{1\leq k,l\leq n}c^{kl}\psi_ku_l\\&
-\frac{h'\psi }{v(1-\psi^2)}|D'u|^2\log |D'u|^2.
\end{split}
\end{align}
Since at $x_0$, $$u_\gamma^2=\psi^2(1+|Du|^2)=\psi^2(1+|D'u|^2+u_\gamma^2),$$
then
\begin{align}\label{5ugamma^2}
 u_\gamma^2=\frac{\psi^2}{1-\psi^2}(1+|D'u|^2).
\end{align}
If
\begin{align}\label{5a0|D'u|^2}
a_0|D'u|^2<u_\gamma^2, \quad a_0=\max_{x\in\partial\Omega}\frac{2\psi^2}{1-\psi^2},
\end{align}
then we get the estimates
\begin{align}\label{5a0|D'u|^21}
(a_0\frac{1-\psi^2}{\psi^2}-1) |D'u|^2<1,
\quad |D'u|^2<\frac{1}{a_0\frac{1-\psi^2}{\psi^2}-1},
\end{align}
and we complete this proof.

So we can assume
\begin{align}\label{5a0|D'u|^22}
a_0|D'u|^2\geq u_\gamma^2,
\end{align}
 then from  $|Du|^2=|D'u|^2+u_\gamma^2$, we have
\begin{align}\label{5a0|D'u|^23}
 |Du|^2\leq (1+a_0)|D'u|^2.
 \end{align}
Now we assume at $x_0$,  we have
\begin{align}\label{5a0|D'u|^24}
|Du|\geq \max\{10\sqrt{(1+a_0)}, 2\sqrt{n}\max_{x\in\partial\Omega}{\frac{|\psi|}{1-\psi^2}}\},
 \end{align}
 then we can get the the following estimates at $x_0$,
 \begin{align}\label{5a0|D'u|^25}
 |D'u|\geq \max\{10, \frac{2\sqrt{n}}{\sqrt{1+a_0}}\max_{x\in\partial\Omega}{\frac{|\psi|}{1-\psi^2}}\}.
 \end{align}

 Inserting \eqref{5a0|D'u|^25} into \eqref{5varphigamma1},  by the choice of $h(u), g(d)$ in \eqref{5hg1},  it follows that at $x_0$,
\begin{align}\label{5varphigamma2}
\begin{split}
|D'u|^2\log|D'u|^2\frac{\partial\varphi}{\partial\gamma}
\geq &\big[\tilde{\alpha}_0-\frac{\sqrt{n}\tilde{\alpha}_0}{v}\frac{|\psi|}{1-\psi^2}
-2\sqrt{1+a_0}\frac{|\psi|_{C^1(\partial\Omega)}}{\sqrt{1-\psi^2}}-C_0\big]|D'u|^2\log|D'u|^2\\
\geq &|D'u|^2\log|D'u|^2\\
>&0.
\end{split}
\end{align}
 On the other hand, by the Hopf Lemma, we have

  $$\frac{\partial\varphi}{\partial\gamma}(x_0)\leq 0,$$
  it is a contradiction to \eqref{5varphigamma2}.
Then we have
$$|D'u|(x_0)\leq\max\{10, \max_{x\in\partial\Omega}{\frac{1}{\sqrt{a_0\frac{1-\psi^2}{\psi^2}-1}}}, \frac{2\sqrt{n}}{\sqrt{1+a_0}}\max_{x\in\partial\Omega}{\frac{|\psi|}{1-\psi^2}}\}.$$

{\bf Case 2.} $ x_0\in \partial\Omega_{\mu_{0}}\bigcap\Omega$. This is due to interior gradient estimates. From Remark~\ref{Rem1.1}, we have
 \begin{align}\label{5CaseII}
\sup_{\partial\Omega_{\mu_0}\bigcap\Omega}|Du|\leq \tilde{M}_1.
\end{align}
where $\tilde{M}_1$ is a positive constant depending only on $n, M_0, \mu_0, L_1$.

{\bf Case 3.} $x_0\in\Omega_{\mu_{0}}$. \par
As in the proof of the Case III in Theorem~\ref{Thm1.1},   we can let $0<\mu_{0}<\mu_1$ be sufficiently small positive constant. As before we have
$$
\sup_{\Omega}|Du|^2\leq C_1(1+\sup_{\partial\Omega}|Du|^2),$$
 where $C_1$ is a positive constant depending on $n, L_1, M_0$.
From  Case 1,   we can assume \eqref{5a0|D'u|^22}, otherwise we have finished the proof of  Theorem~\ref{Thm1.2}.
Hence,
\begin{align}\label{5sup|Du|^22}
\begin{split}
\sup_{\Omega}|Du|^2\leq& C_1\big[1+(1+a_0)\sup_{\partial\Omega\bigcap\{|D'u|\geq 1\}}|D'u|^2\big]\\
\leq& C_2\sup_{\partial\Omega\bigcap\{|D'u|\geq 1\}}|D'u|^2.
\end{split}
\end{align}
 So we have
\begin{align}\label{5sup|Du|^23}
\begin{split}
\sup_{\Omega_{\mu_0}}|Du|^2
\leq& C_2\sup_{\overline{\Omega}_{\mu_0}(M)}|D'u|^2,
\end{split}
\end{align}
  where $\Omega_{\mu_0}(M)=\Omega_{\mu_0}\bigcap\{|D'u|\geq M\}$, $M>10$ is a positive constant;
   $C_2$ is a positive constant depending on $n, L_1, b_0, M_0$.

Assume  $x_1\in\Omega_{\mu_0}(M)$ such that
\begin{align}\label{5sup|D'u|^21}
\begin{split}
\sup_{\Omega_{\mu_0}(M)}|D'u|^2
=& |D'u|^2(x_1).
\end{split}
\end{align}
Since $\tilde{P}(x_0)\geq \tilde{P}(x_1)$, then we have
\begin{align}\label{5sup|D'u|^22}
\begin{split}
\log|D'u|^2(x_0)h(u(x_0))g(d(x_0))\geq&\log|D'u|^2(x_1)h(u(x_1))g(d(x_1)).
\end{split}
\end{align}
It follows that
\begin{align}\label{5sup|D'u|^23}
\begin{split}
|D'u|^2(x_1)\leq&C_3|D'u|^2(x_0).
\end{split}
\end{align}
where  $C_3$ is a positive constant depending on $n, L_1, b_0, L_2, M_0$.
However
\begin{align}\label{5sup|Du|^24}
\begin{split}
\sup_{\Omega_{\mu_0}}|Du|^2
\leq& C_2\sup_{\Omega_{\mu_0}(M)}|D'u|^2=|D'u|^2(x_1)\leq C_3|D'u|^2(x_0).
\end{split}
\end{align}
 Assume $|D'u|(x_0)\geq M>10$, otherwise we get the estimate. Hence at $x_0$,
\begin{align}\label{5sup|Du|^25}
\begin{split}
u_1^2(x_0)=|Du|^2(x_0)
\leq& C_3c^{11}u_1^2(x_0).
\end{split}
\end{align}
Then we have at $x_0$,
\begin{align}\label{5c^{11}}
\begin{split}
c^{11}
\geq& \frac{1}{C_3}>0.
\end{split}
\end{align}
Similar calculations to Case III in the proof of Theorem~\ref{Thm1.1}., and by  the choice of $h(u), g(d)$ in \eqref{5hg1}, we can obtain at last
\begin{align}\label{5aijvarphiij}
0\geq &\sum_{1\leq i,j\leq n}a^{ij}\varphi_{ij}\notag\\
\geq & \big\{h'^2+[(c^{11})^2(n-2)-c^{11}(n-3)]g'^2\big\}u_1^2\log u_1-C_{4}u_1^2\notag\\
\geq & 3u_1^2\log u_1-C_{5}u_1^2.
\end{align}
There exists a positive constant
$C_{6}$ such that
 \begin{align}\label{5C6}
 |D'u|(x_0)\leq C_{6}.
 \end{align}
So from Case 1, Case 2, Case 3, we complete the new proof of Theorem~\ref{Thm1.2}.\qed\\

\begin{Rem}\label{Rem5.1}
 For $q>1$ or $q=0$,   $\psi=\psi(x,u)$,   under the condition $|\psi_u|\leq O(v^{q-1}\log v)$,
 we can  obtain the gradient estimate for the problem \eqref{1.1}-\eqref{1.2}.
 \end{Rem}

\section{Proof of  Theorem~\ref{Thm1.3}.}
{\em Proof of  Theorem~\ref{Thm1.3}.}

In order to unify the computation with the proof of Theorem~\ref{Thm1.1}, we still use the summation index from $1$ to $n$, and at last we take $n=2,3$.

As in the proof of Theorem~\ref{Thm1.1}, let

 $$\Phi(x)=\log|D'u|^2e^{1+M_0+u}e^{\beta_0 d}, \quad  x \in \overline \Omega_{\mu_{0}},$$ where $\beta_0=2L_3+C_0+1 $, $C_0$ is  a positive constant depending only on $n,\Omega$.

Set
\begin{align}\label{6phi}
 \phi(x)= \log \Phi(x) =\log\log|D'u|^2+h(u)+g(d).
 \end{align}
and in the Neumann boundary value ,we choose
\begin{align}\label{6hg1}
 h(u)=1+M_0+u, \quad g(d)=\beta_0 d.
 \end{align}

 We assume that
$\phi(x)$ attains its maximum at $x_0 \in \overline \Omega_{\mu_{0}}$, where $0<\mu_0<\mu_1$ is a sufficiently small number which we shall decide it  later.

{\bf Case i.} If $x_0\in \partial \Omega$, we shall get the bound of $|D'u|(x_0)$.

Following the similar calculations to case I in the proof of Theorem~\ref{Thm1.1}., we let $q=1$ in \eqref{3.10}, and get
\begin{align}\label{6phigamma1}
\begin{split}
|D'u|^2\log |D'u|^2 \frac{\partial\phi}{\partial\gamma}(x_0)
=&\big(g'-h'\psi\big)|D'u|^2\log |D'u|^2-2\sum_{1\leq i,k,l\leq n}c^{kl}u_iu_l(\gamma^i)_k\\&-2\sum_{1\leq k,l\leq n}c^{kl}D_k\psi u_l.
\end{split}
\end{align}
Assume $|Du|(x_0)\ge \sqrt{20}$, otherwise we get the estimates. At $x_0$,  since
$$|Du|^2=|D'u|^2+u_\gamma^2=|D'u|^2+|\psi|^2,$$ so we can assume
\begin{align}\label{6|D'u|^2(x_0)}
\begin{split}
|D'u|^2(x_0) \geq |\psi|^2_{C^0(\partial\Omega\times[-M_0, M_0])},
\end{split}
\end{align}
then $$\max\{20,2|\psi|^2_{C^0(\partial\Omega\times[-M_0,\, M_0])}\}\leq|Du|^2(x_0)\leq2|D'u|^2(x_0),$$
otherwise we get the estimates.
So by the above formulas and \eqref{6phigamma1}, we have
\begin{align}\label{6phigamma2}
\begin{split}
|D'u|^2\log |D'u|^2 \frac{\partial\phi}{\partial\gamma}(x_0)
\geq &\big(\beta_0-|\psi|-|\nabla_x\psi|-2|\psi_u|-C_0\big)|D'u|^2\log |D'u|^2\\
\geq &\big(\beta_0-2L_3-C_0\big)|D'u|^2\log |D'u|^2\\
\geq &|D'u|^2\log |D'u|^2\\
>& 0.
\end{split}
\end{align}
On the other hand, from Hopf Lemma,
$$\frac{\partial\phi}{\partial\gamma}(x_0)\leq 0,$$
it is a contradiction to \eqref{6phigamma2}.

Then we have
\begin{align}\label{6|D'u|(x_0)a}
|D'u|(x_0)\leq\max\{\sqrt{10}, |\psi|_{C^0(\partial\Omega\times[-M_0, M_0])}\}.
\end{align}

{\bf Case ii.} $ x_0\in \partial\Omega_{\mu_{0}}\bigcap\Omega$. This is due to interior gradient estimates. From Remark~\ref{Rem1.1}, we have
 \begin{align}\label{6Case ii}
\sup_{\partial\Omega_{\mu_0}\bigcap\Omega}|Du|\leq \tilde{M}_1.
\end{align}
where $\tilde{M}_1$ is a positive constant depending only on $n, M_0, \mu_0, L_1$.\par

{\bf Case iii.} $x_0\in\Omega_{\mu_{0}}$. \par
As in the proof of the Case III in Theorem ~\ref{Thm1.1},   we can let $0<\mu_{0}<\mu_1$ be sufficiently small positive constant. As  we have
$$
\sup_{\Omega}|Du|^2\leq C_1(1+\sup_{\partial\Omega}|Du|^2),$$
 where $C_1$ is a positive constant depending on $n, L_1, M_0$.
From  Case i,   we can assume \eqref{6|D'u|^2(x_0)}, otherwise we have finished the proof of  Theorem~\ref{Thm1.3}.
Hence,
\begin{align}\label{6sup|Du|^22}
\sup_{\Omega}|Du|^2\leq& C_1\big[1+2\sup_{\partial\Omega\bigcap\{|D'u|\geq 1\}}|D'u|^2\big]
\leq C_2\sup_{\partial\Omega\bigcap\{|D'u|\geq 1\}}|D'u|^2.
\end{align}
So as before, there exists a positive constant $C_2$ such that at $x_0$
\begin{align}\label{6c^{11}}
\begin{split}
c^{11}
\geq& \frac{1}{C_2}>0.
\end{split}
\end{align}
By the above choice of $h, g $ in \eqref{6hg1} and from the similar calculations to the proof of Theorem ~\ref{Thm1.1}, we have for $n=2,3$,
\begin{align}\label{6aijphiij}
\begin{split}
0\geq &\sum_{1\leq i,j\leq n}a^{ij}\phi_{ij}(x_0)\\
\geq & \big\{h'^2+[(c^{11})^2(n-2)-c^{11}(n-3)]\beta_0^2\big\}u_1^2\log u_1-C_{3}u_1^2\\
\geq & u_1^2\log u_1-C_{3}u_1^2.
\end{split}
\end{align}
So there exists  $C_{4}$ such that  $$|D'u|(x_0)\leq C_{4}.$$
Where  the above positive constants $C_{2}, C_2, C_4$ are  depending only on $n, \Omega, \mu_0, M_0, L_1, L_3$.
As in the proof of Theorem~\ref{Thm1.1}, combining three cases, we finally get the following estimate
$$\sup_{\overline\Omega_{\mu_0}}|Du|\leq \max\{M_1, M_2\},$$
where positive constant $M_1$ depends only on $n, \mu_0, M_0, L_1$; $M_2$ depends only on $n, \Omega, \mu_0, M_0,$\\
$ L_1, L_3$.\par
So we complete the new proof of Theorem~\ref{Thm1.3}.\qed\\

{\bf Acknowledgement:} The paper is one part of the author's thesis in University of Science and Technology of China, she would like  to thank Professor Ma Xi-Nan for his guides.

\end{document}